\documentclass{article}

\usepackage{graphicx}
\usepackage{float}
\usepackage{amsthm}
\usepackage{amsmath}
\usepackage{amssymb}
\usepackage{mathrsfs}
\usepackage{dsfont}
\usepackage{tikz-cd}
\usepackage{stmaryrd}
\usepackage{enumitem}
\usepackage{placeins}
\usepackage{pdflscape}
\usepackage{pdflscape}
\usepackage{hyperref}

\theoremstyle{definition}
\newtheorem{Definition}{Definition}[subsection]

\theoremstyle{plain}
\newtheorem{Theorem}[Definition]{Theorem}

\theoremstyle{plain}
\newtheorem{Proposition}[Definition]{Proposition}

\theoremstyle{plain}
\newtheorem{Lemma}[Definition]{Lemma}

\theoremstyle{plain}
\newtheorem{Corollary}[Definition]{Corollary}

\theoremstyle{plain}

\theoremstyle{definition}
\newtheorem{Example}[Definition]{Example}

\theoremstyle{remark}
\newtheorem{Remark}[Definition]{Remark}

\theoremstyle{plain}
\newcommand{\thistheoremname}{}
\newtheorem*{genericthm*}{\thistheoremname}
\newenvironment{namedthm*}[1]
  {\renewcommand{\thistheoremname}{#1}%
   \begin{genericthm*}}
  {\end{genericthm*}}

\title{Local Modules in Braided Monoidal 2-Categories}
\author{Thibault D. D\'{e}coppet and Hao Xu}
\date{July 2023}

\begin{document}

\bibliographystyle{alpha}

\maketitle

\begin{abstract}
Given an algebra in a monoidal 2-category, one can construct a 2-category of right modules. Given a braided algebra in a braided monoidal 2-category, it is possible to refine the notion of right module to that of a local module. Under mild assumptions, we prove that the 2-category of local modules admits a braided monoidal structure. In addition, if the braided monoidal 2-category has duals, we go on to show that the 2-category of local modules also has duals. Furthermore, if it is a braided fusion 2-category, we establish that the 2-category of local modules is a braided multifusion 2-category. We examine various examples. For instance, working within the 2-category of 2-vector spaces, we find that the notion of local module recovers that of braided module 1-category. Finally, we examine the concept of a Lagrangian algebra, that is a braided algebra with trivial 2-category of local modules. In particular, we completely describe Lagrangian algebras in the Drinfeld centers of fusion 2-categories, and we discuss how this result is related to the classifications of topological boundaries of (3+1)d topological phases of matter.
\end{abstract}

\section*{Introduction}

It is well-know that the 1-category of modules over a commutative algebra is symmeric monoidal. The notion of a module over an algebra can be internalized to any monoidal 1-category $\mathcal{C}$. Further, provided we work in a braided monoidal 1-category $\mathcal{B}$, it is also sensible to consider a commutative algebra $B$ in $\mathcal{B}$. Under mild assumptions on $\mathcal{B}$, the 1-category of $B$-modules in $\mathcal{B}$, which we denote by $\mathbf{Mod}_{\mathcal{B}}(B)$, admits a canonical monoidal structure given by the relative tensor product over $B$. However, $\mathbf{Mod}_{\mathcal{B}}(B)$ is not braided unless $\mathcal{B}$ is symmetric. It was nevertheless observed in \cite{P} that the full sub-1-category $\mathbf{Mod}_{\mathcal{B}}^{loc}(B)$ of $\mathbf{Mod}_{\mathcal{B}}(B)$ on the local modules, also known as dyslectic modules, admits a braiding.

One noteworthy application of local modules comes from its relation with Drinfeld centers \cite{Sch}. More precisely, given a commutative algebra $B$ in $\mathcal{Z}(\mathcal{C})$, the Drinfeld center of the monoidal 1-category $\mathcal{C}$, it is shown under mild assumptions on $\mathcal{C}$ that $\mathbf{Mod}_{\mathcal{Z}(\mathcal{C})}^{loc}(B)$ is equivalent as a braided monoidal 1-category to $\mathcal{Z}(\mathbf{Mod}_{\mathcal{C}}(B))$.
Local modules were also used in \cite{KO} as a way to produce new modular tensor 1-categories from the known examples. We note that this last problem was initially undertaken using subfactors \cite{BEK99}. Subsequently, the relation between 1-categories of local modules and non-degenerate braided fusion 1-categories was explored much further in \cite{DMNO}. For instance, it was established that if $\mathcal{B}$ is a non-degenerate braided fusion 1-category, and $B$ is a connected separable commutative algebra, also called a connected étale algebra, in $\mathcal{B}$, then $\mathbf{Mod}^{loc}_{\mathcal{B}}(B)$ is a non-degenerate braided fusion 1-category. These results were then generalized in \cite{DNO}.

Categorifying the notion of a fusion 1-category, fusion 2-categories were introduced in \cite{DR}. The theory of algebras in fusion 2-categories was then extensively developed in \cite{D4, D7, D8, D9} by the first author. Further, it was established in \cite{DY} that the 2-category $\mathbf{Mod}_{\mathfrak{B}}(B)$ of modules over a braided separable algebra $B$ in a braided fusion 2-category $\mathfrak{B}$ is a multifusion 2-category. In this context, braided separable algebras played a central role in the proof of the minimal non-degenerate extension conjecture \cite{JFR}.

Motivations for developing the theory of local modules in braided fusion 1-categories also come from condensed matter Physics. Specifically, the theory of local modules plays a crucial role in the physical theory of anyon condensations, a program which was initiated in \cite{BSS02,BSS03} and was fully realized in \cite{K14}. More precisely, a modular tensor 1-category $\mathcal{B}$ describes a (2+1)-dimensional topological order (up to invertible ones), a connected étale algebra $B$ in $\mathcal{B}$ represents a combination of anyons that condense to the vacuum in a new phase. The modular tensor 1-category associated to the new phase, also called the 1-category of deconfined particles, is precisely $\mathbf{Mod}_\mathcal{B}^{loc}(B)$. A more detailed account of the historical development of anyon condensation theory can be found in \cite{K14}.

Going up in dimension, local modules in braided fusion 2-categories are expected to play a role in the study of (3+1)-dimensional anyon condensation. Some work has already been dome in this direction, specifically in the (3+1)-dimensional toric code model \cite{ZLZHKT}, where examples of Lagrangian algebras are used to describe its topological boundaries. %We will come back to the classification of \'etale algebras and Lagrangian algebras in braided fusion 2-categories in our future work.

\subsection*{Results}

Our first objective is to categorify the main result of \cite{P}. We fix a braided monoidal 2-category $\mathfrak{B}$ and a braided algebra $B$ in $\mathfrak{B}$. We review the definition of local right $B$-modules in $\mathfrak{B}$ introduced in \cite{ZLZHKT} using the variant of the graphical calculus of \cite{GS} introduced in \cite{D4}. More precisely, a local right $B$-module is a right $B$-module equipped with a 2-isomorphism, called a holonomy. Relying on \cite{D8} and \cite{DY}, we establish the following result.

\begin{namedthm*}{Theorem \ref{thm:2catlocalmodulesbraided}}
Let $\mathfrak{B}$ be a braided monoidal 2-category, and $B$ a braided algebra in $\mathfrak{B}$. If $\mathfrak{B}$ has relative tensor products over $B$ and they are preserved by the monoidal product of $\mathfrak{B}$, then the 2-category $\mathbf{Mod}^{loc}_{\mathfrak{B}}(B)$ of local $B$-modules admits a braided monoidal structure.
\end{namedthm*}

\noindent We note that a similar result was obtained independently in \cite{Pom} using multi-2-categories.

We go on to study more specifically the case when $\mathfrak{B}$ is a braided multifusion 2-category. Under these hypotheses, the relative tensor product exists provided we consider separable algebras as introduced in \cite{JFR} (see also \cite{D7} for a thorough discussion). Thus, we take $B$ to be an étale algebra in $\mathfrak{B}$, that is a braided separable algebra in $\mathfrak{B}$. We show that the underlying 2-category of $\mathbf{Mod}^{loc}_{\mathfrak{B}}(B)$ is finite semisimple. Further, we prove that the dual of a right $B$-module admits a compatible holonomy. Putting these facts together yield the next theorem, which categorifies a number of results of \cite{KO}.

\begin{namedthm*}{Theorem \ref{thm:localmodulesmultifusion}}
Let $B$ be an étale algebra in a braided multifusion 2-category $\mathfrak{B}$. Then, $\mathbf{Mod}_{\mathfrak{B}}^{loc}(B)$ is a braided multifusion 2-category.
\end{namedthm*}

\noindent Further, if $\mathfrak{B}$ is a fusion 2-category, and $B$ is a connected étale algebra, we find that $\mathbf{Mod}_{\mathfrak{B}}^{loc}(B)$ is a braided fusion 2-category.

Then, we identify the braided fusion 2-category $\mathbf{Mod}_{\mathfrak{B}}^{loc}(B)$ when $\mathfrak{B}$ is a braided fusion 2-category of interest. Firstly, when $\mathfrak{B}=\mathbf{2Vect}$, connected étale algebras are precisely braided fusion 1-categories. Given a braided fusion 1-category $\mathcal{B}$, we show that $\mathbf{Mod}_{\mathbf{2Vect}}^{loc}(\mathcal{B})$ is equivalent as a braided monoidal 2-category to $\mathscr{Z}(\mathbf{Mod}(\mathcal{B}))$, the Drinfeld center of $\mathcal{B}$. We do so by proving that local right $\mathcal{B}$-modules in $\mathbf{2Vect}$ correspond exactly to the finite semisimple braided $\mathcal{B}$-module 1-categories of \cite{DN}. Secondly, we prove that taking local modules is invariant under base change, i.e.\ given a 1-morphism between \'etale algebras $A \rightarrow B$ in a braided fusion 2-category $\mathfrak{B}$, one has that the 2-category of local $B$-modules in $\mathfrak{B}$ is equivalent to the 2-category of local $B$-modules in $\mathbf{Mod}_{\mathfrak{B}}^{loc}(A)$.

We use our previous results to study Lagrangian algebras, that is connected étale algebras whose associated 2-category of local modules is trivial. Conceptually, this last condition should be thought of as a categorical non-degeneracy condition. Namely, Lagrangian algebras in $\mathbf{2Vect}$ are exactly given by non-degenerate braided fusion 1-categories. We compare our notion of Lagrangian algebra with that introduced in \cite{JFR}, and show that, given any braided fusion 1-category $\mathcal{B}$, Lagrangian algebras in $\mathscr{Z}(\mathbf{Mod}(\mathcal{B}))$ correspond exactly to non-degenerate braided fusion 1-categories equipped with a braided functor from $\mathcal{B}$. Finally, we end by discussing the relation between our abstract work on Lagrangian algebras, and the classification of topological boundaries of (3+1)d topological phases of matter. In the (3+1)d toric code model, we compare our results with those obtained from lattice model considerations in \cite{ZLZHKT}.

\subsection*{Acknowledgements}

We would like to thank Liang Kong and Matthew Yu for discussions, as well as feedback on a draft of this manuscript. H.X. was supported by DAAD Graduate School Scholarship Programme (57572629) and DFG Project 398436923.

\section{Preliminaries}

\subsection{Graphical Calculus}

We work within a monoidal 2-category $\mathfrak{C}$ with monoidal product $\Box$ and monoidal unit $I$ in the sense of \cite{SP}. Thanks to the coherence theorem of \cite{Gur}, we may assume without loss of generality that $\mathfrak{C}$ is strict cubical in the sense of definition 2.26 of \cite{SP}. We use the graphical calculus of \cite{GS}, as described in \cite{D4} and \cite{D7}. We often omit the symbol $\Box$ from our notations, and we use the symbol $1$ to denote an identity 1-morphism. The interchanger is depicted using the string diagram below on the left, and its inverse by that on the right: $$\begin{tabular}{c c c c}
\includegraphics[width=20mm]{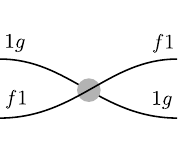},\ \ \ \   & \ \ \ \  \includegraphics[width=20mm]{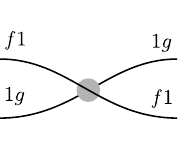}.
\end{tabular}$$ In particular, the lines correspond to 1-morphisms, and the coupons to 2-morphisms. The regions represent objects, which are uniquely determined by the 1-morphisms. Further, our string diagrams are read from top to bottom, which yields the compositions of 1-morphisms, and then from left to right.

For our purposes, it is also necessary to recall the graphical conventions related to 2-natural transformations from \cite{GS}. These will only be used for the braiding, which will be introduced below. Let $F,G:\mathfrak{A}\rightarrow \mathfrak{B}$ be two (weak) 2-functors, and let $\tau:F\Rightarrow G$ be 2-natural transformation. That is, for every object $A$ in $\mathfrak{A}$, we have a 1-morphism $\tau_A:F(A)\rightarrow G(A)$, and for every 1-morphism $f:A\rightarrow B$ in $\mathfrak{A}$, we have a 2-isomorphism $$\begin{tikzcd}[sep=tiny]
F(A) \arrow[ddd, "F(f)"']\arrow[rrr, "\tau_A"]  &                                        &    & G(A) \arrow[ddd, "G(f)"]  \\
 &  &    & \\
  &  &  &  \\
F(B)\arrow[rrr, "\tau_B"']\arrow[Rightarrow, rrruuu, "\tau_f", shorten > = 2ex, shorten < = 2ex]                                            &                                        &    &  G(B), 
\end{tikzcd}$$
These 2-isomorphisms have to satisfy obvious coherence relations. In our graphical language, we will depict the 2-isomorphism $\tau_f$ using the following diagram on the left, and its inverse using the diagram on the right: $$\begin{tabular}{c c c c}
\includegraphics[width=20mm]{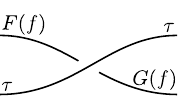},\ \ \ \   & \ \ \ \  \includegraphics[width=20mm]{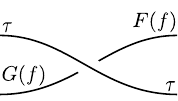}.
\end{tabular}$$ 

\subsection{Braided Monoidal 2-Categories}

In the present article, we will for the most part work within $\mathfrak{B}$ a braided monoidal 2-category in the sense of \cite{SP}. Thanks to the coherence theorem of \cite{Gur2}, we may assume that $\mathfrak{B}$ is a semi-strict braided monoidal 2-category. In particular, $\mathfrak{B}$ comes equipped with a braiding $b$, which is an adjoint 2-natural equivalence given on objects $A,B$ in $\mathfrak{B}$ by $$b_{A,B}:A\Box B\rightarrow B\Box A.$$ Its pseudo-inverse will be denoted by $b^{\bullet}$. Further, there are two invertible modifications $R$ and $S$, which are given on the objects $A,B,C$ of $\mathfrak{B}$ by
\begin{center}
\begin{tabular}{@{}c c@{}}

$\begin{tikzcd}
ABC \arrow[rr, "b"] \arrow[rd, "b1"'] & {} \arrow[d, Rightarrow, "R"]          & BCA, \\
                                      & BAC \arrow[ru, "1b"'] &    
\end{tikzcd}$

&

$\begin{tikzcd}
ABC \arrow[rr, "b_2"] \arrow[rd, "1b"'] & {} \arrow[d, Rightarrow, "S"]     & CAB \\
                                        & ACB \arrow[ru, "b1"'] &    
\end{tikzcd}$
\end{tabular}
\end{center}

\noindent where the subscript in $b_2$ records were the braiding occurs. To avoid any possible confusion, we will systematically write $b$ instead of a would be $b_1$ as this can too easily be confused with $b1$. Further, these modifications are subject to the following relations, which are taken from section 2.1.1 of \cite{DY}:

\begin{enumerate}
\item [a.] For every objects $A,B,C,D$ in $\mathfrak{B}$, we have
\end{enumerate}

\newlength{\calculus}

\settoheight{\calculus}{\includegraphics[width=37.5mm]{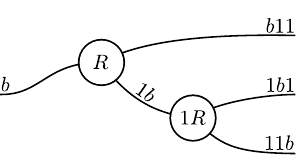}}

\begin{equation}\label{eqn:braidingaxiom1}
\begin{tabular}{@{}ccc@{}}

\includegraphics[width=37.5mm]{Pictures/Preliminaries/Braiding/braidingaxiom1.pdf} & \raisebox{0.45\calculus}{$=$} &
\includegraphics[width=37.5mm]{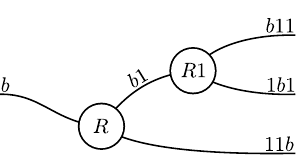}

\end{tabular}
\end{equation}

\begin{enumerate}
\item [] in $Hom_{\mathfrak{B}}(ABCD, BCDA)$,
\item [b.] For every objects $A,B,C,D$ in $\mathfrak{B}$, we have
\end{enumerate}

\settoheight{\calculus}{\includegraphics[width=37.5mm]{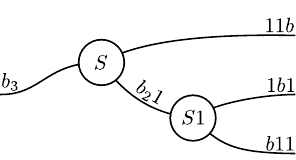}}

\begin{equation}\label{eqn:braidingaxiom2}
\begin{tabular}{@{}ccc@{}}

\includegraphics[width=37.5mm]{Pictures/Preliminaries/Braiding/braidingaxiom3.pdf} & \raisebox{0.45\calculus}{$=$} &
\includegraphics[width=37.5mm]{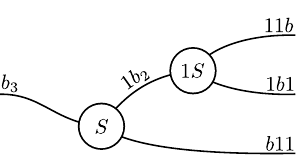}

\end{tabular}
\end{equation}

\begin{enumerate}
\item [] in $Hom_{\mathfrak{B}}(ABCD, DABC)$,
\item [c.] For every objects $A,B,C,D$ in $\mathfrak{B}$, we have
\end{enumerate}

\settoheight{\calculus}{\includegraphics[width=37.5mm]{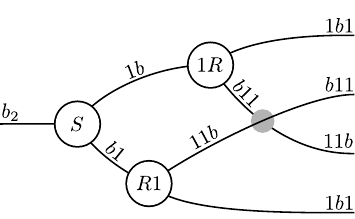}}

\begin{equation}\label{eqn:braidingaxiom3}
\begin{tabular}{@{}ccc@{}}

\includegraphics[width=45mm]{Pictures/Preliminaries/Braiding/braidingaxiom5.pdf} & \raisebox{0.45\calculus}{$=$} &
\includegraphics[width=37.5mm]{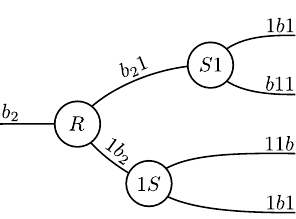}

\end{tabular}
\end{equation}

\begin{enumerate}
\item [] in $Hom_{\mathfrak{B}}(ABCD, CDAB)$,
\item [d.] For every objects $A,B,C$ in $\mathfrak{B}$, we have
\end{enumerate}

\settoheight{\calculus}{\includegraphics[width=45mm]{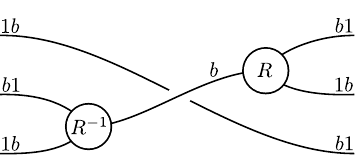}}

\begin{equation}\label{eqn:braidingaxiom4}
\begin{tabular}{@{}ccc@{}}

\includegraphics[width=45mm]{Pictures/Preliminaries/Braiding/braidingaxiom7.pdf} & \raisebox{0.45\calculus}{$=$} &
\includegraphics[width=45mm]{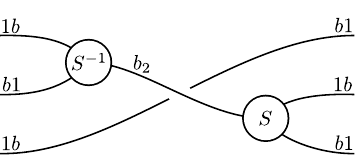}

\end{tabular}
\end{equation}

\begin{enumerate}
\item [] in $Hom_{\mathfrak{B}}(ABC, CBA)$,
\item [e.] For every object $A$ in $\mathfrak{B}$, the adjoint 2-natural equivalences $$b_{A,I}:A\Box I\rightarrow I\Box A\textrm{ and } b_{I,A}:I\Box A\rightarrow A\Box I$$ are the identity adjoint 2-natural equivalences,

\item [f.] For every objects $A,B,C$ in $\mathfrak{B}$, the 2-isomorphisms $R_{A,B,C}$ and $S_{A,B,C}$ are the identity 2-isomorphism whenever either $A$, $B$, or $C$ is equal to $I$.
\end{enumerate}

Let us now examine some examples. For simplicity, we work over an algebraically closed field $\mathds{k}$ of characteristic zero, but this is not necessary \cite{D5}.

\begin{Example}
We write $\mathbf{2Vect}$ for the 2-category of finite semisimple 1-categories. The Deligne tensor product endows $\mathbf{2Vect}$ with a symmetric monoidal structure. This is the most fundamental example of a symmetric fusion 2-category as introduced by \cite{DR} (see also \cite{D2} for a slightly different perspective). More generally, we can also consider $\mathbf{FinCat}$, the 2-category of finite 1-categories, and right exact functors. The Deligne tensor product also endows $\mathbf{FinCat}$ with a symmetric monoidal structure.
\end{Example}

\begin{Example}
Let $\mathcal{E}$ be a symmetric fusion 1-category over $\mathds{k}$. Recall from \cite{D5} that the 2-category $\mathbf{Mod}(\mathcal{E})$ of finite semisimple right $\mathcal{E}$-module 1-categories is a fusion 2-category with monoidal structure given by the relative Deligne tensor product. Further, it admits a symmetric monoidal structure \cite{DY}.
\end{Example}

\begin{Example}
Let $G$ be a finite group. We can consider the symmetric fusion 2-category $\mathbf{2Rep}(G)$ of finite semisimple 1-categories with a $G$-action. In fact, it was shown in \cite{D8} that $\mathbf{2Rep}(G)\simeq\mathbf{Mod}(\mathbf{Rep}(G))$ as symmetric fusion 2-categories. More generally, one can consider the 2-category of 2-representations of a finite 2-group.
\end{Example}

\begin{Example}
Let $\mathfrak{C}$ be a fusion 2-category. Then, the Drinfeld center of $\mathfrak{C}$, which we denote by $\mathscr{Z}(\mathfrak{C})$, is a braided monoidal 2-category \cite{Cr}. Further, it was shown in \cite{D9} that $\mathscr{Z}(\mathfrak{C})$ is a fusion 2-category. For instance, given a finite group $G$, one can consider the braided fusion 2-category $\mathscr{Z}(\mathbf{2Vect}_G)$, which was described explicitly in \cite{KTZ}. We will be particularly interested in the case when $\mathfrak{C}$ is connected, i.e.\ $\mathfrak{C}\simeq \mathbf{Mod}(\mathcal{B})$ for some braided fusion 1-category $\mathcal{B}$. Thanks to the results of \cite{D9}, this is not a loss of generality as the Drinfeld center of any fusion 2-category is of this form.
\end{Example}

\subsection{Algebras and Modules}

Let $\mathfrak{C}$ be a strict cubical monoidal 2-category. We recall the definition of an algebra in $\mathfrak{C}$ from \cite{D7}. The definition of an algebra in an arbitrary monoidal 2-category using our graphical conventions may be found in \cite{D4}.

\begin{Definition}\label{def:algebra}
An algebra in $\mathfrak{C}$ consists of:
\begin{enumerate}
    \item An object $A$ of $\mathfrak{C}$;
    \item Two 1-morphisms $m:A\Box A\rightarrow A$ and $i:I\rightarrow A$;
    \item Three 2-isomorphisms
\end{enumerate}
\begin{center}
\begin{tabular}{@{}c c c@{}}
$\begin{tikzcd}[sep=small]
A \arrow[rrrr, equal] \arrow[rrdd, "i1"'] &  & {} \arrow[dd, Rightarrow, "\lambda"', near start, shorten > = 1ex] &  & A \\
                                   &  &                           &  &   \\
                                   &  & AA, \arrow[rruu, "m"']     &  &  
\end{tikzcd}$

&

$\begin{tikzcd}[sep=small]
AAA \arrow[dd, "1m"'] \arrow[rr, "m1"]    &  & AA \arrow[dd, "m"] \\
                                            &  &                      \\
AA \arrow[rr, "m"'] \arrow[rruu, Rightarrow, "\mu", shorten > = 2.5ex, shorten < = 2.5ex] &  & A,                   
\end{tikzcd}$

&

$\begin{tikzcd}[sep=small]
                                  &  & AA \arrow[rrdd, "m"] \arrow[dd, Rightarrow, "\rho", shorten > = 1ex, shorten < = 2ex] &  &   \\
                                  &  &                                             &  &   \\
A \arrow[rruu, "1i"] \arrow[rrrr,equal] &  & {}                                          &  & A,
\end{tikzcd}$

\end{tabular}
\end{center}

satisfying:

\begin{enumerate}
\item [a.] We have:
\end{enumerate}

\newlength{\prelim}

\settoheight{\prelim}{\includegraphics[width=52.5mm]{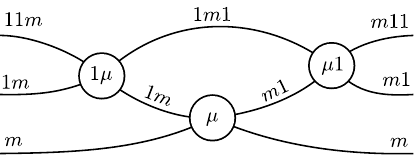}}

\begin{equation}\label{eqn:algebraassociativity}
\begin{tabular}{@{}ccc@{}}

\includegraphics[width=52.5mm]{Pictures/Preliminaries/Algebra/associativity1.pdf} & \raisebox{0.45\prelim}{$=$} &
\includegraphics[width=40mm]{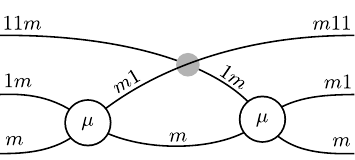},

\end{tabular}
\end{equation}

\begin{enumerate}
\item [b.] We have:
\end{enumerate}

\settoheight{\prelim}{\includegraphics[width=22.5mm]{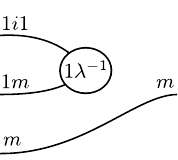}}

\begin{equation}\label{eqn:algebraunitality}
\begin{tabular}{@{}ccc@{}}

\includegraphics[width=22.5mm]{Pictures/Preliminaries/Algebra/unitality1.pdf} & \raisebox{0.45\prelim}{$=$} &

\includegraphics[width=37.5mm]{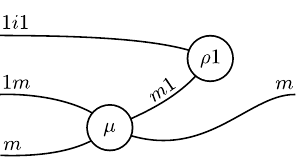}.

\end{tabular}
\end{equation}
\end{Definition}

Let us now recall the definition of a right $A$-module in $\mathfrak{C}$ from definition 1.2.3 of \cite{D7}. We invite the reader to consult definition 3.2.1 of \cite{D4} for the definition in a general monoidal 2-category.

\begin{Definition}\label{def:module}
A right $A$-module in $\mathfrak{C}$ consists of:
\begin{enumerate}
    \item An object $M$ of $\mathfrak{C}$;
    \item A 1-morphism $n^M:M\Box A\rightarrow M$;
    \item Two 2-isomorphisms
\end{enumerate}
\begin{center}
\begin{tabular}{@{}c c@{}}
$\begin{tikzcd}[sep=small]
MAA \arrow[dd, "1m"'] \arrow[rr, "n^M1"]    &  & MA \arrow[dd, "n^M"] \\
                                            &  &                      \\
MA \arrow[rr, "n^M"'] \arrow[rruu, Rightarrow, "\nu^M", shorten > = 2.5ex, shorten < = 2.5ex] &  & M,                   
\end{tikzcd}$

&

$\begin{tikzcd}[sep=small]
                                  &  & MA \arrow[rrdd, "n^M"] \arrow[dd, Rightarrow, "\rho^M", shorten > = 1ex, shorten < = 2ex] &  &   \\
                                  &  &                                             &  &   \\
M \arrow[rruu, "1i"] \arrow[rrrr,equal] &  & {}                                          &  & M,
\end{tikzcd}$
\end{tabular}
\end{center}

satisfying:

\begin{enumerate}
\item [a.] We have:
\end{enumerate}

\settoheight{\prelim}{\includegraphics[width=52.5mm]{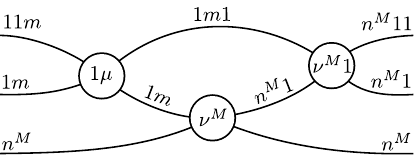}}

\begin{equation}\label{eqn:moduleassociativity}
\begin{tabular}{@{}ccc@{}}

\includegraphics[width=52.5mm]{Pictures/Preliminaries/Module/associativity1.pdf} & \raisebox{0.45\prelim}{$=$} &
\includegraphics[width=45mm]{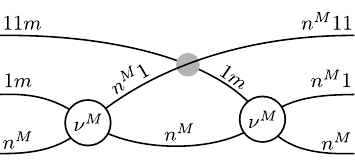},

\end{tabular}
\end{equation}

\begin{enumerate}
\item [b.] We have:
\end{enumerate}

\settoheight{\prelim}{\includegraphics[width=22.5mm]{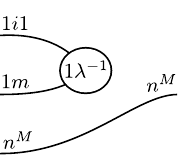}}

\begin{equation}\label{eqn:moduleunitality}
\begin{tabular}{@{}ccc@{}}

\includegraphics[width=22.5mm]{Pictures/Preliminaries/Module/unitality1.pdf} & \raisebox{0.45\prelim}{$=$} &

\includegraphics[width=37.5mm]{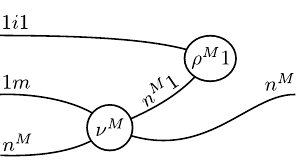}.

\end{tabular}
\end{equation}
\end{Definition}

Finally, let us recall definitions 3.2.6 and 3.2.7 of \cite{D4}.

\begin{Definition}\label{def:modulemap}
Let $M$ and $N$ be two right $A$-modules. A right $A$-module 1-morphism consists of a 1-morphism $f:M\rightarrow N$ in $\mathfrak{C}$ together with an invertible 2-morphism

$$\begin{tikzcd}[sep=small]
MA \arrow[dd, "f1"'] \arrow[rr, "n^M"]    &  & M \arrow[dd, "f"] \\
                                            &  &                      \\
NA \arrow[rr, "n^N"'] \arrow[rruu, Rightarrow, "\psi^f", shorten > = 2.5ex, shorten < = 2.5ex] &  & N,                   
\end{tikzcd}$$

subject to the coherence relations:

\begin{enumerate}
\item [a.] We have:
\end{enumerate}

\settoheight{\prelim}{\includegraphics[width=52.5mm]{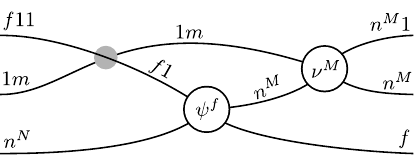}}

\begin{equation}\label{eqn:modulemapassociativity}
\begin{tabular}{@{}ccc@{}}

\includegraphics[width=52.5mm]{Pictures/Preliminaries/Module/map1.pdf} & \raisebox{0.45\prelim}{$=$} &

\includegraphics[width=52.5mm]{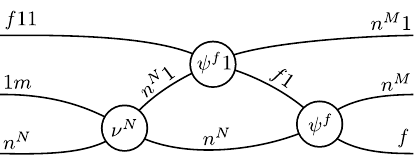},

\end{tabular}
\end{equation}

\begin{enumerate}
\item [b.] We have:
\end{enumerate}

\settoheight{\prelim}{\includegraphics[width=30mm]{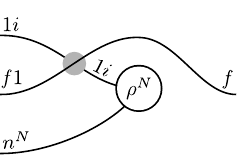}}

\begin{equation}\label{eqn:modulemapunitality}
\begin{tabular}{@{}ccc@{}}

\includegraphics[width=30mm]{Pictures/Preliminaries/Module/map3.pdf} & \raisebox{0.45\prelim}{$=$} &

\includegraphics[width=30mm]{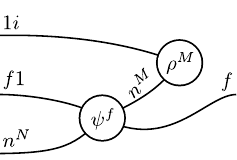}.

\end{tabular}
\end{equation}
\end{Definition}

\begin{Definition}\label{def:moduleintertwiner}
Let $M$ and $N$ be two right $A$-modules, and $f,g:M\rightarrow M$ two right $A$-module 1-morphisms. A right $A$-module 2-morphism $f\Rightarrow g$ is a 2-morphism $\gamma:f\Rightarrow g$ in $\mathfrak{C}$ that satisfies the following equality:

\settoheight{\prelim}{\includegraphics[width=30mm]{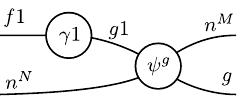}}

$$\label{eqn:module2map}\begin{tabular}{@{}ccc@{}}

\includegraphics[width=30mm]{Pictures/Preliminaries/Module/2morphism1.pdf} & \raisebox{0.45\prelim}{$=$} &

\includegraphics[width=30mm]{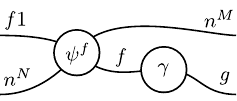}.

\end{tabular}$$
\end{Definition}

These object assemble to form a 2-category as was proven in \cite{D4}.

\begin{Lemma}
Let $A$ be an algebra in a monoidal 2-category $\mathfrak{C}$. Right $A$-modules, right $A$-module 1-morphisms, and right $A$-module 2-morphisms form a 2-category, which we denote by $\mathbf{Mod}_{\mathfrak{C}}(A)$.
\end{Lemma}

\begin{Example}
In $\mathbf{2Vect}$, algebras correspond exactly to finite semisimple monoidal 1-categories. Given a finite semisimple monoidal 1-category $\mathcal{C}$, right $\mathcal{C}$-modules in $\mathbf{2Vect}$ are precisely finite semisimple right $\mathcal{C}$-module 1-categories. A similar observation holds for module morphisms. More generally, algebras $\mathbf{FinCat}$ are precisely finite monoidal 1-categories, whose monoidal product is right exact in both variables. Fixing such a monoidal 1-category $\mathcal{C}$, right $\mathcal{C}$-modules in $\mathbf{FinCat}$ correspond exactly to finite right $\mathcal{C}$-module 1-categories, for which the action is right exact in both variables.
\end{Example}

Let us also recall the following definition from \cite{D4}.

\begin{Definition}\label{def:algebrahomomorphism}
Let $A$ and $B$ be two algebras in $\mathfrak{C}$. An algebra 1-homomorphism $f:A\rightarrow B$ consists of a 1-morphism $f:A\rightarrow B$ in $\mathfrak{C}$, together with two invertible 2-morphisms
 
\begin{center}
\begin{tabular}{@{}c c@{}}
$$\begin{tikzcd}[sep=tiny]
AA \arrow[rrr, "f1"] \arrow[dddd, "m^A"'] & & & BA \arrow[rrr, "1f"] &                                     &    & BB \arrow[dddd, "m^B"] \\
                                                  &&&    &  &      &                             \\
                                                   &       {} \arrow[rrrr, "\kappa^f", Rightarrow, shorten <=4ex, shorten >=4ex]         &&&&                   {}    &                              \\
                                                   &               &&&                     &  &                              \\
A \arrow[rrrrrr, "f"']                                &         &&  &                          &    & B,                         
\end{tikzcd}$$

&

$$\begin{tikzcd}[sep=tiny]
I \arrow[rrr, equal] \arrow[ddd, "i^A"'] &                                     &    & I \arrow[ddd, "i_B"] \\
                                                     &  &  {}   &                              \\
                                                   &  {}                            &  &                              \\
A \arrow[rrr, "f"'] \arrow[uuurrr, "\iota^f", Rightarrow, shorten <=4ex, shorten >=4ex]                               &                                     &    & B,                         
\end{tikzcd}$$
\end{tabular}
\end{center}

satisfying:

\begin{enumerate}
\item [a.] We have:
\end{enumerate}

\settoheight{\prelim}{\includegraphics[width=45mm]{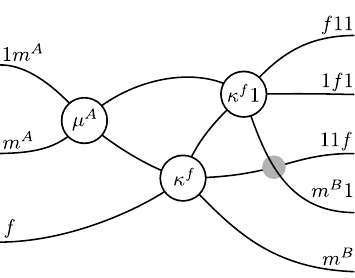}}

\begin{equation}\label{eqn:algebrahomomorphismassociativity}
\begin{tabular}{@{}ccc@{}}

\includegraphics[width=45mm]{Pictures/Preliminaries/Algebra/morphismassociativity1.pdf} & \raisebox{0.45\prelim}{$=$} &

\includegraphics[width=45mm]{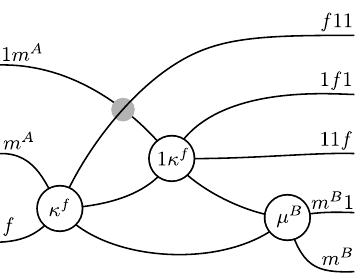},

\end{tabular}
\end{equation}

\begin{enumerate}
\item [b.] We have:
\end{enumerate}

\settoheight{\prelim}{\includegraphics[width=45mm]{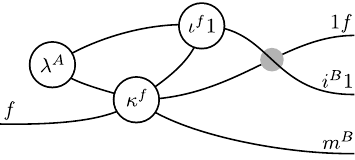}}

\begin{equation}\label{eqn:algebrahomomorphismunitelity1}
\begin{tabular}{@{}ccc@{}}

\includegraphics[width=45mm]{Pictures/Preliminaries/Algebra/morphismunitality1.pdf} & \raisebox{0.45\prelim}{$=$} &

\includegraphics[width=30mm]{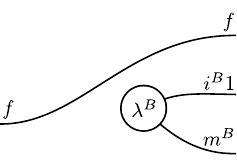},

\end{tabular}
\end{equation}

\begin{enumerate}
\item [c.] We have:
\end{enumerate}

\settoheight{\prelim}{\includegraphics[width=45mm]{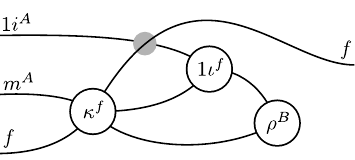}}

\begin{equation}\label{eqn:algebrahomomorphismunitality2}
\begin{tabular}{@{}ccc@{}}

\includegraphics[width=45mm]{Pictures/Preliminaries/Algebra/morphismunitality3.pdf} & \raisebox{0.45\prelim}{$=$} &

\includegraphics[width=30mm]{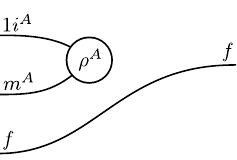}.

\end{tabular}
\end{equation}
\end{Definition}

\subsection{Relative Tensor Product}\label{sub:tensor}

Recall that we are working with a strict cubical monoidal 2-category $\mathfrak{C}$. Let us fix an algebra $A$ in $\mathfrak{C}$. In addition, let $M$ be a right $A$-module in $\mathfrak{C}$, and $N$ be a left $A$-module $N$ in $\mathfrak{C}$ (for which we use the notations of \cite{D7}). Following \cite{D8}, we can examine whether the relative tensor product of $M$ and $N$ over $A$ exists. This is determined by a 2-universal property, which we now recall for later use.

\begin{Definition}\label{def:balanced1morphism}
Let $C$ be an object of $\mathfrak{C}$. An $A$-balanced 1-morphism $(M,N)\rightarrow C$ consists of:
\begin{enumerate}
    \item A 1-morphism $f:M\Box N\rightarrow C$ in $\mathfrak{C}$;
    \item A 2-isomorphism
\end{enumerate}

$$\begin{tikzcd}[sep=small]
MAN \arrow[dd, "1l^N"'] \arrow[rr, "n^M1"]    &  & MN \arrow[dd, "f"] \\
                                            &  &                      \\
MN \arrow[rr, "f"'] \arrow[rruu, Rightarrow, "\tau^f", shorten > = 2.5ex, shorten < = 2.5ex] &  & C,                   
\end{tikzcd}$$

satisfying:

\begin{enumerate}
\item [a.]
\end{enumerate}

\newlength{\tensor}
\settoheight{\tensor}{\includegraphics[width=52.5mm]{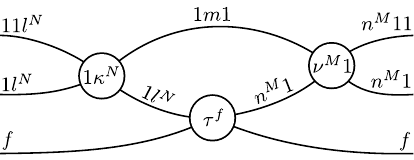}}

\begin{equation}\label{eqn:balancedassociativity}
\begin{tabular}{@{}ccc@{}}

\includegraphics[width=52.5mm]{Pictures/Preliminaries/Tensor/associativity1.pdf} & \raisebox{0.45\tensor}{$=$} &
\includegraphics[width=45mm]{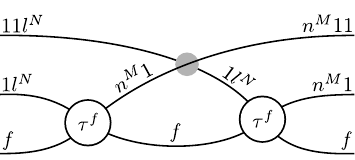},

\end{tabular}
\end{equation}

\begin{enumerate}
\item [b.]
\end{enumerate}

\settoheight{\tensor}{\includegraphics[width=22.5mm]{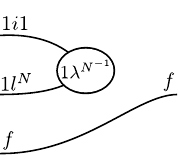}}

\begin{equation}\label{eqn:balancedunitality}
\begin{tabular}{@{}ccc@{}}

\includegraphics[width=22.5mm]{Pictures/Preliminaries/Tensor/unitality1.pdf} & \raisebox{0.45\tensor}{$=$} &

\includegraphics[width=37.5mm]{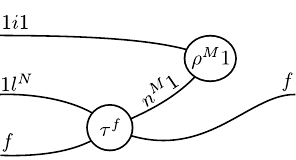}.

\end{tabular}
\end{equation}
\end{Definition}

\begin{Definition}\label{def:balanced2morphism}
Let $C$ be an object of $\mathfrak{C}$, and $f,g:(M,N)\rightarrow C$ be two $A$-balanced 1-morphisms. An $A$-balanced 2-morphism $f\Rightarrow g$ is a 2-morphism $\gamma:f\Rightarrow g$ in $\mathfrak{C}$ such that
\settoheight{\tensor}{\includegraphics[width=30mm]{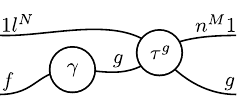}}

\begin{center}
\begin{tabular}{@{}ccc@{}}

\includegraphics[width=30mm]{Pictures/Preliminaries/Tensor/2morphism1.pdf} & \raisebox{0.45\tensor}{$=$} &

\includegraphics[width=30mm]{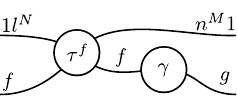}.

\end{tabular}
\end{center}
\end{Definition}

\begin{Definition}\label{def:relativetensor}
The relative tensor product of $M$ and $N$ over $A$, if it exists, is an object $M\Box_A N$ of $\mathfrak{C}$ together with an $A$-balanced 1-morphism $t_A:(M,N)\rightarrow M\Box_A N$ satisfying the following 2-universal property:
\begin{enumerate}
    \item For every $A$-balanced 1-morphism $f:(M,N)\rightarrow C$, there exists a 1-morphism $\widetilde{f}:M\Box_A N\rightarrow C$ in $\mathfrak{C}$ and an $A$-balanced 2-isomorphism $\xi:\widetilde{f}\circ t_A\cong f$.
    \item For any 1-morphisms $g,h:M\Box_A N\rightarrow C$ in $\mathfrak{C}$, and any $A$-balanced 2-morphism $\gamma:g\circ t_A\Rightarrow h\circ t_A$, there exists a unique 2-morphism $\zeta:g\Rightarrow h$ such that $\zeta\circ t_A = \gamma$.
\end{enumerate}
\end{Definition}

\begin{Example}\label{ex:relativetensor2Vect}
Let us fix an algebra in $\mathbf{2Vect}$, that is a finite semisimple monoidal 1-category $\mathcal{C}$. Recall that we are working over an algebraically closed field $\mathds{k}$ of characteristic zero for simplicity. Further, let us fix a finite semisimple right $\mathcal{C}$-module 1-category $\mathcal{M}$, and a finite semisimple left $\mathcal{C}$-module 1-category $\mathcal{N}$. We view $\mathcal{M}$ as a right $\mathcal{C}$-module in $\mathbf{2Vect}$, and $\mathcal{N}$ as a left $\mathcal{C}$-module. In this case, the relative tensor product $\mathcal{M}\boxtimes_{\mathcal{C}}\mathcal{N}$ is precisely the relative Deligne tensor product of \cite{ENO}. Namely, the above definitions correspond exactly to the definitions of section 3.1 of \cite{ENO}. The aforementioned reference shows that the relative tensor product exists if we assume that $\mathcal{C}$ has duals, i.e.\ $\mathcal{C}$ is a multifusion 1-category. For completeness, let us also note that this example can be generalized to $\mathbf{FinCat}$ (see \cite{DSPS14}).

For later use, it is convenient to compare in more detail the explicit construction of the relative Deligne tensor product given in section 3.2 of \cite{DN} with the abstract categorical definition above. More precisely, an object of $\mathcal{M}\boxtimes_{\mathcal{C}}\mathcal{N}$ is a pair $(V,\gamma)$ consisting of an object $V$ of $\mathcal{M}\boxtimes\mathcal{N}$ together with natural isomorphisms $$\gamma_C: V\otimes(C\boxtimes I)\rightarrow V\otimes(I\boxtimes C),$$ for every $C$ in $\mathcal{C}$ satisfying the obvious coherence relations. Then, the canonical $\mathcal{C}$-balanced functor $t:\mathcal{M}\boxtimes\mathcal{N}\rightarrow \mathcal{M}\boxtimes_{\mathcal{C}}\mathcal{N}$ is the right adjoint of the forgetful functor. Moreover, given objects $M$ in $\mathcal{M}$ and $N$ in $\mathcal{N}$, and writing $t(M\boxtimes N) = (V,\gamma)$, the balancing natural isomorphism $\tau^t$ on $M\boxtimes N$ is given by $$V\otimes(I\boxtimes C)\xrightarrow{\gamma_C^{-1}} V\otimes(I\boxtimes C).$$ It follows from the definitions that this construction satisfies the desired 2-universal property.
\end{Example}

Generalizing the above example, it is interesting to ask when the relative tensor product of $M$ and $N$ over $A$ exists. We will be particularly interested in the case when $\mathfrak{C}$ is a multifusion 2-category in the sense of \cite{D2}. Under this hypothesis, it was shown in \cite{D8} that a sufficient condition for the existence of the relative tensor product is that the algebra $A$ be separable, a notion which has its origin in \cite{GJF}. This recovers the previous example as separable algebra in the fusion 2-category $\mathbf{2Vect}$ are exactly multifusion 1-categories. For completeness, we also recall the intermediate definition of a rigid algebra due to \cite{G}. Both of these notions were studied extensively in \cite{D7}, to which we refer the reader for further discussion and examples.

\begin{Definition}
A rigid algebra in a monoidal 2-category is an algebra $A$ such that the 1-morphism $m:A\Box A\rightarrow A$ admits a right adjoint $m^*$ as an $A$-$A$-bimodule 1-morphism.
\end{Definition}

\begin{Definition}
A separable algebra in a monoidal 2-category is a rigid algebra $A$ such that the counit $\epsilon^m:m\circ m^*\Rightarrow Id_A$ splits as an $A$-$A$-bimodule 2-morphism.
\end{Definition}

Provided that the relative tensor product over $A$ of any two modules exists in $\mathfrak{C}$ and commutes with the monoidal structure, the 2-category $\mathbf{Bimod}_{\mathfrak{C}}(A)$ of $A$-$A$-bimodules in $\mathfrak{C}$ inherits a monoidal structure given by $\Box_A$, the relative tensor product over $A$ \cite{D8}.
In particular, this is the case if $A$ is a separable algebra in a fusion 2-category $\mathfrak{C}$. In this case, the 2-category $\mathbf{Bimod}_{\mathfrak{C}}(A)$ of $A$-$A$-bimodules in $\mathfrak{C}$ has a monoidal structure provided by $\Box_A$. In particular, we note that all of the coherence data is supplied by the 2-universal property of the relative tensor product.

\subsection{Braided and \'Etale Algebras}

Let $\mathfrak{B}$ be a semi-strict braided monoidal 2-category. We recall the following definition from \cite{DY}.

\begin{Definition}\label{def:braidedalgebra}
A braided algebra in $\mathfrak{B}$ consists of:
\begin{enumerate}
    \item An algebra $B$ in $\mathfrak{B}$;
    \item A 2-isomorphisms
\end{enumerate}

$$\begin{tikzcd}[sep=small]
                                  &  & BB \arrow[rrdd, "m"] \arrow[dd, Rightarrow, "\beta", shorten > = 1ex, shorten < = 2ex] &  &   \\
                                  &  &                                             &  &   \\
BB \arrow[rruu, "b"] \arrow[rrrr, "m"'] &  & {}                                          &  & B,
\end{tikzcd}$$

satisfying:

\begin{enumerate}
\item [a.] We have:
\end{enumerate}

\newlength{\braid}

\settoheight{\braid}{\includegraphics[width=52.5mm]{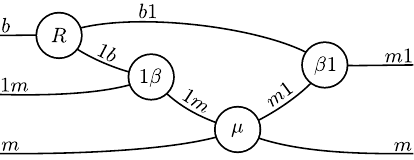}}

\begin{equation}\label{eqn:braidedalgebra1}
\begin{tabular}{@{}ccc@{}}

\includegraphics[width=52.5mm]{Pictures/Preliminaries/BraidedAlgebra/braidedalgebra1.pdf} & \raisebox{0.45\braid}{$=$} &
\includegraphics[width=52.5mm]{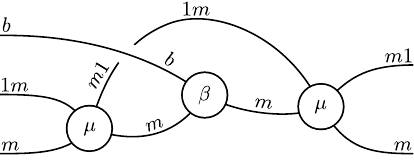},

\end{tabular}
\end{equation}

\begin{enumerate}
\item [b.] We have:
\end{enumerate}

\settoheight{\braid}{\includegraphics[width=52.5mm]{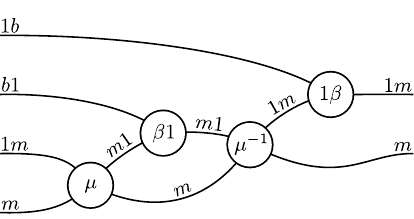}}

\begin{equation}\label{eqn:braidedalgebra2}
\begin{tabular}{@{}ccc@{}}

\includegraphics[width=52.5mm]{Pictures/Preliminaries/BraidedAlgebra/braidedalgebra3.pdf} & \raisebox{0.45\braid}{$=$} &

\includegraphics[width=52.5mm]{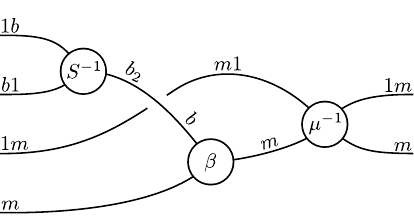},

\end{tabular}
\end{equation}

\begin{enumerate}
\item [c.] We have:
\end{enumerate}

\settoheight{\braid}{\includegraphics[width=30mm]{Pictures/Preliminaries/BraidedAlgebra/braidedalgebra4.pdf}}

\begin{equation}\label{eqn:braidedalgebra3}
\begin{tabular}{@{}ccc@{}}

\includegraphics[width=30mm]{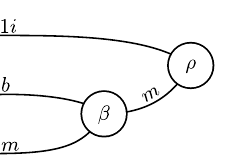} & \raisebox{0.45\braid}{$=$} &

\includegraphics[width=30mm]{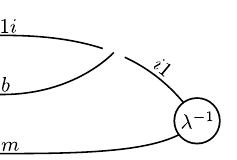}.

\end{tabular}
\end{equation}
\end{Definition}

\begin{Definition}\label{def:braidedalgebrahomomorphism}
Let $A$ and $B$ be two braided algebras in $\mathfrak{B}$. A braided algebra 1-homomorphism $f:A\rightarrow B$ is an algebra 1-homomorphism $f:A\rightarrow B$ that satisfies:

\settoheight{\braid}{\includegraphics[width=45mm]{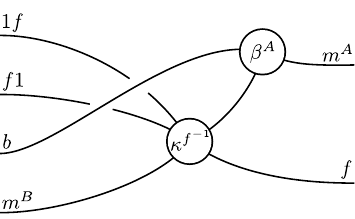}}

$$\begin{tabular}{@{}ccc@{}}

\includegraphics[width=45mm]{Pictures/Preliminaries/BraidedAlgebra/braidedhomomorphism1.pdf} & \raisebox{0.45\braid}{$=$} &

\includegraphics[width=37.5mm]{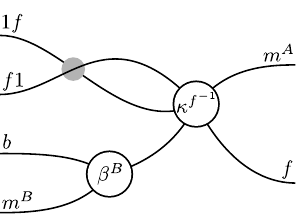}.

\end{tabular}$$
\end{Definition}

Let us now assume that $\mathfrak{B}$ is a braided fusion 2-category over algebraically closed of characteristic zero. We note that the assumptions on the ground field can be relaxed \cite{D5,D7}. We will be particularly interested in the following objects.

\begin{Definition}
An \'etale algebra in a braided multifusion 2-category $\mathfrak{B}$ is a separable braided algebra. An étale algebra $B$ is called connected provided that its unit 1-morphism $i:I\rightarrow B$ is simple, i.e.\ $End_{\mathfrak{B}}(i)\cong \mathds{k}$.
\end{Definition}

\begin{Example}
Braided rigid algebras in $\mathbf{FinCat}$ are exactly finite braided multitensor 1-categories in the sense of \cite{EGNO}. \'Etale algebras in $\mathbf{2Vect}$ are given by separable braided multifusion 1-categories, and connected étale algebras are braided fusion 1-categories.
\end{Example}

\begin{Example}
Let $\mathcal{E}$ be a symmetric fusion 1-category. It follows from \cite{D9} that étale algebras in $\mathbf{Mod}(\mathcal{E})$ are exactly braided multifusion 1-categories equipped with a symmetric functor from $\mathcal{E}$ to their symmetric center.
\end{Example}

\begin{Example}
Let $G$ be a finite group. Braided rigid algebras in $\mathscr{Z}(\mathbf{2Vect}_G)$ are exactly $G$-crossed braided multifusion 1-categories. Thanks to corollary 5.1.2 of \cite{D9}, all such algebras are étale. We note that a $G$-crossed braided multifusion 1-category $\mathcal{B}$ yields a connected étale algebra in $\mathscr{Z}(\mathbf{2Vect}_G)$ if and only if the canonical $G$-action on the monoidal unit of $\mathcal{B}$ permutes all the simple summands transitively.
\end{Example}

\begin{Example}
Let $\mathcal{B}$ be a braided fusion 1-category. It follows from \cite{DN} and \cite{D9} that étale algebras in $\mathscr{Z}(\mathbf{Mod}(\mathcal{B}))$ are braided multifusion 1-categories equipped with a braided functor from $\mathcal{B}$. Further, connected étale algebras are braided fusion 1-categories equipped with a braided functor from $\mathcal{B}$. Such objects played a key role in the proof of the minimal non-degenerate extension conjecture obtained in \cite{JFR}.
\end{Example}

\subsection{Induction 2-Functors}\label{sub:induction}

Let us now fix a braided algebra $B$ in a braided monoidal 2-category $\mathfrak{B}$, which we assume is semi-strict without loss of generality. Further, we will assume that relative tensor products over $B$ exist in $\mathfrak{B}$, and are preserved by the monoidal product of $\mathfrak{B}$. Under these hypotheses, and building upon \cite{D8}, it was shown in section 3.2 of \cite{DY} that the 2-category of right $B$-modules in $\mathfrak{B}$ admits a monoidal structure. We denote this monoidal 2-category by $\mathbf{Mod}^+_{\mathfrak{B}}(B)$, and its monoidal product by $\Box_B^+$. For later use, let us recall more precisely how this monoidal structure is constructed. Given a right $B$-module $M$, we can define a $B$-$B$-bimodule $Ind^+(M)$ by endowing the right $B$-module $M$ with a compatible left $B$-module structure using the left action $l^M_+:=n^M\circ b:B\Box M\rightarrow M$, and the following 2-isomorphisms 

\settoheight{\braid}{\includegraphics[width=30mm]{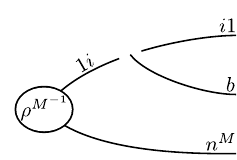}}

$$\raisebox{15pt}{\raisebox{0.45\braid}{$\lambda^M:=\ $}
\includegraphics[width=30mm]{Pictures/Preliminaries/Induction/lambdaM.pdf},}\ \ \ \settoheight{\braid}{\includegraphics[width=52.5mm]{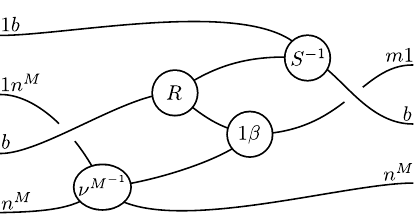}}
\raisebox{0.45\braid}{$\kappa^M:=\ $}
\includegraphics[width=52.5mm]{Pictures/Preliminaries/Induction/kappaM.pdf},$$

\settoheight{\braid}{\includegraphics[width=52.5mm]{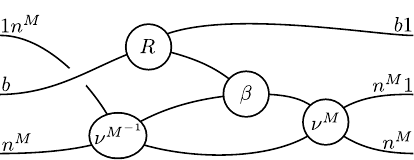}}

$$\raisebox{0.45\braid}{$\beta^M:=\ $}
\includegraphics[width=52.5mm]{Pictures/Preliminaries/Induction/betaM.pdf}.$$

\noindent Here, we are using the notations of \cite{D8} for left modules, and bimodules. Then, it was shown in lemma 3.3 of \cite{DY} that the above assignment extends to a 2-functor $$Ind^+:\mathbf{Mod}^+_{\mathfrak{B}}(B)\rightarrow \mathbf{Bimod}_{\mathfrak{B}}(B),$$ which is fully faithful on 2-morphisms. In particular, we may view $\mathbf{Mod}^+_{\mathfrak{B}}(B)$ as a sub-2-category of $\mathbf{Bimod}_{\mathfrak{B}}(B)$. But, the relative tensor product over $B$ endows the 2-category $\mathbf{Bimod}_{\mathfrak{B}}(B)$ with a monoidal structure \cite{D8}. Further, it was established in proposition 3.4 of \cite{DY} that the sub-2-category $\mathbf{Mod}^+_{\mathfrak{B}}(B)$ is closed under the relative tensor product over $B$. In particular, it admits a monoidal structure, with monoidal product denoted using $\Box_B^+$. For our current purposes, it is convenient to rephrase this result in the following way.

\begin{Lemma}
The 2-functor $Ind^+:\mathbf{Mod}^+_{\mathfrak{B}}(B)\rightarrow \mathbf{Bimod}_{\mathfrak{B}}(B)$ is monoidal.
\end{Lemma}

Given $M$, $N$ two right $B$-modules, we will need to have an explicit description of the equivalence \begin{equation}\label{eqn:Ind+monoidal}X_{M,N}^+:Ind^+(M)\Box_BInd^+(N)\rightarrow Ind^+(M\Box_B^+N),\end{equation} witnessing that $Ind^+$ is monoidal. It follows from the proof of proposition 3.4 of \cite{DY} that the underlying right $B$-module 1-morphism is the identity on $M\Box_B^+N$. Further, the left $B$-module structure is constructed using the 2-universal property of the relative tensor product over $B$. More precisely, if $t:M\Box N\rightarrow M\Box^+_B N$ denotes the 2-universal $B$-balanced right $B$-module 1-morphism with balancing $\tau^t_+:t\circ (M\Box l^N_+)\Rightarrow t\circ (n^M\Box N)$, then the 2-isomorphism $$(t\circ 1n^N \circ R^{-1})\cdot(\tau^{t^{-1}}_+\circ b1):t\circ (n^M\Box N) \circ (b\Box N)\cong t\circ (M\Box n^N)\circ b$$ is $B$-balanced. Thus, there exists a 2-isomorphism $\chi^{X_{M,N}^+}:l^{M\Box^+_BN}\cong n^{M\Box^+_BN}\circ b$ ensuring that $X_{M,N}^+$ is compatible with the left $B$-module structures.

As indicated in \cite{DY}, the above construction admits a variant that uses $b^{\bullet}$ instead of $b$. Succinctly, given a right $B$-module $M$, we can define a $B$-$B$-bimodule $Ind^-(M)$ by endowing the right $B$-module $M$ with a compatible left $B$-module structure using the left action $l^M_-:=n^M\circ b^{\bullet}:B\Box M\rightarrow M$, and appropriate 2-isomorphisms. This assignment extends to a 2-functor $$Ind^-:\mathbf{Mod}^-_{\mathfrak{B}}(B)\rightarrow \mathbf{Bimod}_{\mathfrak{B}}(B).$$ Further, there is a monoidal structure on the 2-category $\mathbf{Mod}^-_{\mathfrak{B}}(B)$ with monoidal product $\Box_B^-$, and such that the 2-functor $Ind^-$ is monoidal. Given $M$, $N$ two right $B$-modules, the equivalence \begin{equation}\label{eqn:Ind-monoidal}X^-_{M,N}:Ind^-(M)\Box_BInd^-(N)\rightarrow Ind^-(M\Box_B^-N)\end{equation} has underlying right $B$-module 1-morphism the identity on $M\Box_B^-N$, and its left $B$-module structure is constructed using the 2-universal property of the relative tensor product over $B$.

\section{Local Modules and Their Properties}

\subsection{Definition}

Let us fix $B$ a braided algebra in a braided monoidal 2-category $\mathfrak{B}$. Without loss of generality, we will assume that $\mathfrak{B}$ is semi-strict. We begin by recalling the definition of a local $B$-module, which has already appeared in appendix B of \cite{ZLZHKT} (see also \cite{Pom} for a related definition in a slightly different context).

\begin{Definition} \label{def:localmodule}
A local $B$-module in $\mathfrak{B}$ consists of:
\begin{enumerate}
    \item A right $B$-module $M$ in $\mathfrak{B}$,
\end{enumerate}

\newlength{\algebra}
\settoheight{\algebra}{\includegraphics[width=75mm]{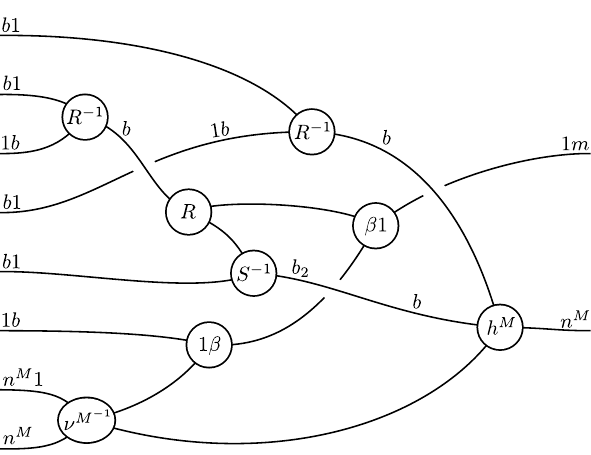}}

\begin{landscape}
        \vspace*{2.2cm}
        \begin{figure}[!htb]
        \begin{equation}\label{eqn:holonomyassociativity}
        \begin{tabular}{@{}ccc@{}}
        \includegraphics[width=75mm]{Pictures/Preliminaries/LocalModule/holonomycoh1left.pdf}&\raisebox{0.48\algebra}{$=$} &
        \includegraphics[width=60mm]{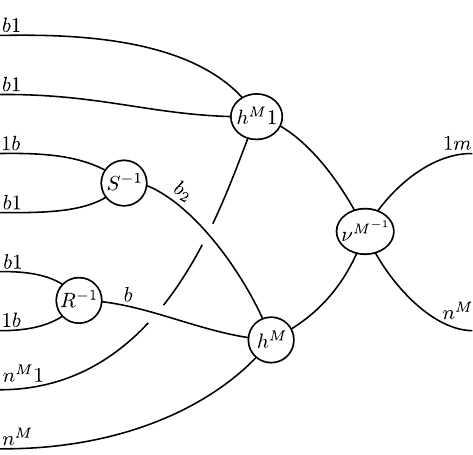},
        \end{tabular}
        \end{equation}
        \end{figure}
\end{landscape}

\begin{enumerate}
    \item[2.] A 2-isomorphism, called a holonomy,
    $$\begin{tikzcd}[sep=small]
    BM \arrow[rr, "b"]                    & {} \arrow[dd, Rightarrow, "h^{M}", shorten > = 2ex, shorten < = 1ex, near start] & MB \arrow[dd, "n^M"] \\
                                      &                        &                      \\
    MB \arrow[uu, "b"] \arrow[rr, "n^M"'] & {}                     & M,                  
    \end{tikzcd}$$
\end{enumerate}

    satisfying:
    \begin{enumerate}
        \item [a.] We have that equation (\ref{eqn:holonomyassociativity}) is satisfied,
    \end{enumerate}

    \begin{enumerate}
        \item [b.] We have:
     \end{enumerate}   
        \settoheight{\algebra}{\includegraphics[width=60mm]{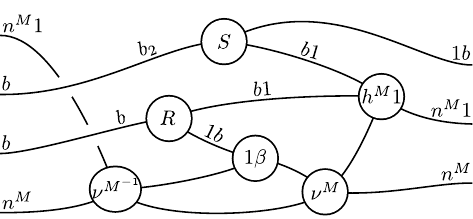}}

        \begin{equation}\label{eqn:holonomyassociativityprime}
        \begin{tabular}{@{}ccc@{}}
        \includegraphics[width=57mm]{Pictures/Preliminaries/LocalModule/holonomycoh1primeleft.pdf}&\raisebox{0.45\algebra}{$=$} &
        \includegraphics[width=50mm]{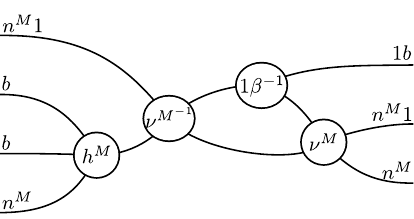},
        \end{tabular}
        \end{equation}
        
    \begin{enumerate}
        \item [c.] We have:
     \end{enumerate}   
        \settoheight{\algebra}{\includegraphics[width=30mm]{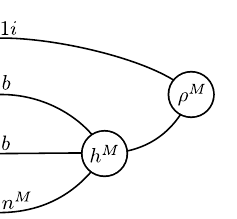}}

        \begin{equation}\label{eqn:holonomyunitality}
        \begin{tabular}{@{}ccc@{}}
        \includegraphics[width=30mm]{Pictures/Preliminaries/LocalModule/holonomycoh2left.pdf}&\raisebox{0.45\algebra}{$=$} &
        \includegraphics[width=30mm]{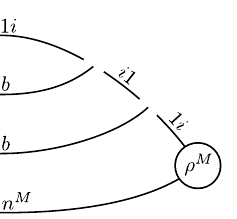}.
        \end{tabular}
        \end{equation}
\end{Definition}

\begin{Definition} \label{def:localmorphism}
Given $M$ and $N$ two local $B$-modules. A 1-morphism of local $B$-modules is a right $B$-module 1-morphism $f$ in $\mathfrak{B}$ satisfying the following equation
    
\settoheight{\algebra}{\includegraphics[width=40mm]{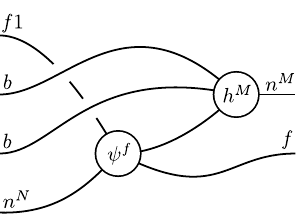}}

\begin{equation}\label{eqn:localmorphism}
    \begin{tabular}{@{}ccc@{}}
        \includegraphics[width=37.5mm]{Pictures/Preliminaries/LocalModule/localmodule1morcohleft.pdf}&\raisebox{0.45\algebra}{$=$} &
        \includegraphics[width=37.5mm]{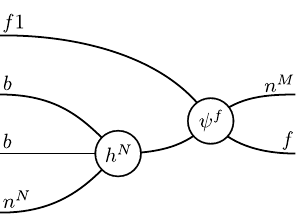}.
    \end{tabular}
\end{equation}
        
\end{Definition}

\begin{Definition} \label{def:local2morphism}
Given $f,g:M\rightarrow N$ two 1-morphisms of local $B$-modules. A 2-morphisms of local $B$-modules is a right $B$-module 2-morphism $f\Rightarrow g$.
\end{Definition}

\begin{Lemma}
Local $B$-modules in $\mathfrak{B}$, 1-morphisms of local $B$-modules, and 2-morphisms of local $B$-modules form a 2-category, which we denote by $\mathbf{Mod}^{loc}_{\mathfrak{B}}(B)$.
\end{Lemma}

The following technical results will play a key role in all of the subsequent constructions. The first lemma is based on an observation given in remark C.8 of \cite{ZLZHKT}. It allows us to give an equivalent perspective on holonomies, which extremely useful in practice.

\begin{Lemma}\label{lem:equivalentperspectiveholonomy}
Let $M$ be a right $B$-module. The data of a holonomy on $M$ corresponds exactly to the data necessary to upgrade the canonical right $B$-module 1-morphism $Id_M:M\rightarrow M$ to a 1-morphism of $B$-$B$-bimodules $\hbar^M:Ind^+(M)\rightarrow Ind^-(M)$.
\end{Lemma}
\begin{proof}
In the notations of section 1.3 of \cite{D7}, the relevant coherence 2-isomorphism witnessing that the $\hbar^M$ is compatible with the left $B$-module structures is of the form $$\begin{tikzcd}[sep=small]
BM \arrow[dd, "b^{\bullet}"'] \arrow[rr, "b"]    &  & MB \arrow[dd, "n^M"] \\
                                            &  &                      \\
MB \arrow[rr, "n^M"'] \arrow[rruu, Rightarrow, "\chi^{\hbar^M}", shorten > = 2.5ex, shorten < = 2.5ex] &  & B,                   
\end{tikzcd}$$ In particular, $\chi^{\hbar^M}$ is obtained from $(h^M)^{-1}$ by currying, i.e.\ we have
\settoheight{\algebra}{\includegraphics[width=30mm]{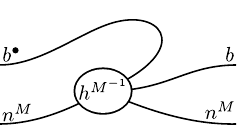}}

$$\begin{tabular}{@{}cc@{}}
    \raisebox{0.45\algebra}{$\chi^{\hbar^M}:=$} &
    \includegraphics[width=37.5mm]{Pictures/Preliminaries/LocalModule/chihM.pdf}.
\end{tabular}$$ It is easy to check that the equations required of $h^M$ to define a holonomy on $M$ correspond exactly to the equations required of $\chi^{\hbar^M}$ to define a compatible left $B$-module structure on $Id_M$ as in the statement of the lemma.
\end{proof}

\begin{Lemma}\label{lem:equivalentperspectiveholonomymorphisms}
Let $M$ and $N$ be two local right $B$-modules, and let $f:M\rightarrow N$ be a right $B$-module 1-morphism. Then, $f$ is a 1-morphism of local $B$-modules if and only if the following square of $B$-$B$-bimodules 1-morphisms commutes $$\begin{tikzcd}[sep=small]
Ind^+(M) \arrow[rr, "\hbar^M"] \arrow[dd, "Ind^+(f)"'] &  & Ind^-(M) \arrow[dd, "Ind^-(f)"] \\
                                            &  &                                 \\
Ind^+(N) \arrow[rr, "\hbar^N"']                         &  & Ind^-(N).                       
\end{tikzcd}$$
\end{Lemma}
\begin{proof}
Provided the square commutes strictly, the proof follows readily by direct inspection. It therefore only remains to explain why it is enough to check that the square commutes weakly. To see this, note that the underlying diagram of right $B$-modules commutes strictly. Moreover, the $B$-$B$-bimodule 1-morphisms $\hbar^M$ and $\hbar^N$ are isomorphisms. Thus, provided that the above square commutes weakly, the choices of filling correspond to the invertible elements in the algebra $End_{B\mathrm{-}B}(Ind^+(f))$. But, this algebra is isomorphic to $End_B(f)$ as the 2-functor $Ind^+$ is fully faithful on 2-morphisms by lemma 3.3. of \cite{DY}. In particular, if the square commutes up to any $B$-$B$-bimodule 2-isomorphism, then it necessarily commutes strictly.
\end{proof}

\subsection{The Braided Monoidal 2-Category of Local Modules}

We show that the 2-category of local modules admits a canonical braided monoidal structure. Our construction proceeds by enhancing the construction of a monoidal 2-category of right modules over a braided algebra given in \cite{DY}. We also point the reader towards \cite{Pom} for related results in a different context. Throughout, we let $\mathfrak{B}$ be a braided monoidal 2-category, and $B$ a braided algebra in $\mathfrak{B}$ such that the relative tensor product over $B$ of any two modules exists. Further, we assume that these relative tensor products are preserved by the monoidal product of $\mathfrak{B}$.

\begin{Proposition} \label{prop:2catlocalmodulemonoidal}
The 2-category $\mathbf{Mod}^{loc}_{\mathfrak{B}}(B)$ of local $B$-modules inherits a monoidal structure from $\mathbf{Mod}^+_{\mathfrak{B}}(B)$.
\end{Proposition}
\begin{proof}
Without loss of generality, we may assume that $\mathfrak{B}$ is semi-strict. As was already recalled in section \ref{sub:induction}, the 2-category $\mathbf{Mod}^+_{\mathfrak{B}}(B)$ of right $B$-modules carries a monoidal structure given by the relative tensor product $\Box^+_B$. Thus, it is enough to enhance this construction to include holonomies. We do so using the perspective of lemma \ref{lem:equivalentperspectiveholonomy}.

More precisely, let $M$ and $N$ be two local right $B$-modules with holonomies $h^M$ and $h^N$ respectively. We will equivalently think of $h^M$ and $h^N$ as $B$-$B$-bimodule 1-morphisms $$\hbar^M:Ind^+(M)\rightarrow Ind^-(M),\ \hbar^N:Ind^+(N)\rightarrow Ind^-(N).$$ Then, the relative tensor product $M\Box^+_BN$ can be equipped with a natural holonomy, namely the one that corresponds to the $B$-$B$-bimodules 1-morphism \begin{align}\label{eqn:holonomyproduct}Ind^+(M\Box_B^+N)&\xrightarrow{(X^+_{M,N})^{\bullet}} Ind^+(M)\Box_BInd^+(N)\xrightarrow{\hbar^M\Box_BId}Ind^-(M)\Box_BInd^+(N)\notag\\
&\xrightarrow{Id\Box_B\hbar^N} Ind^-(M)\Box_BInd^-(N)\xrightarrow{X^-_{M,N}} Ind^-(M\Box^+_BN).\end{align} In order to appeal to lemma \ref{lem:equivalentperspectiveholonomy}, we have to make sure that the underlying right $B$-module 1-morphism is the identity. In order to ensure this, we make the convention that $M\Box_B^+N = M\Box_B^-N$ as right $B$-modules. Namely, the $B$-$B$-bimodules $Ind^+(N)$ and $Ind^-(N)$ are isomorphic via $\hbar^N$. Then, the right $B$-modules $M\Box_B^+N$ equipped with the $B$-balanced right $B$-module 1-morphism $t_{M,N}\circ (M\Box (\hbar^N)^{-1})$ satisfies the 2-universal property of $M\Box_B^-N$. We may therefore without loss of generality assume that they are equal. This convention ensures that the underlying right $B$-module 1-morphism of the composite given in equation (\ref{eqn:holonomyproduct}) is the identity as the underlying right $B$-module of $Ind^-(M)\Box_BInd^-(N)$ is $M\Box_B^-N$ by construction.

Now, let $M$, $N$, and $P$ be right $B$-modules equipped with holonomies $h^M$, $h^N$, and $h^P$ respectively. Thanks to proposition 3.4 of \cite{DY}, there is 2-natural equivalence $$\alpha_{M,N,P}:(M\Box^+_BN)\Box^+_BP\simeq M\Box^+_B(N\Box^+_BP)$$ witnessing the associativity of the monoidal product of $\mathbf{Mod}_{\mathfrak{B}}^+(B)$. In fact, this 2-natural equivalence is compatible with the holonomies, as can be readily seen from lemma \ref{lem:equivalentperspectiveholonomymorphisms} and the fact that $Ind^+$ and $Ind^-$ are monoidal 2-functors. A similar argument shows that the 2-natural equivalences $$l_N:B\Box^+_BN\rightarrow N,\ r_M:M\Box^+_BB\rightarrow M$$ witnessing the unitality of the monoidal product of $\mathbf{Mod}_{\mathfrak{B}}^+(B)$ are also compatible with the holonomies.

Finally, there is nothing to check at the level of invertible modifications. Namely, a 2-morphism of local right $B$-modules is nothing but a 2-morphism of right $B$-modules. They satisfy the necessary equations thanks to proposition 3.4 of \cite{DY}. This finishes the proof of the proposition.
\end{proof}

\begin{Remark}\label{rem:monoidalforgetful}
It follows from the construction of the monoidal structure on $\mathbf{Mod}^{loc}_{\mathfrak{B}}(B)$ that the canonical forgetful 2-functor $$\mathbf{T}:\mathbf{Mod}^{loc}_{\mathfrak{B}}(B)\rightarrow \mathbf{Mod}^{+}_{\mathfrak{B}}(B)$$ admits a monoidal structure. Further, this 2-functor is fully faithful on 2-morphisms. For later use, let us also record that $\mathbf{T}$ is conservative on 1-morphisms, i.e.\ a 1-morphism $f:M\rightarrow N$ is an equivalence in $\mathbf{Mod}^{loc}_{\mathfrak{B}}(B)$ if and only if $\mathbf{T}(f)$ is an equivalence.
\end{Remark}

\begin{Theorem}\label{thm:2catlocalmodulesbraided}
Let $\mathfrak{B}$ be a braided monoidal 2-category, and $B$ a braided algebra in $\mathfrak{B}$. If $\mathfrak{B}$ has relative tensor products over $B$, and they are preserved by the monoidal product of $\mathfrak{B}$, then the 2-category $\mathbf{Mod}^{loc}_{\mathfrak{B}}(B)$ of local $B$-modules admits a braided monoidal structure.
\end{Theorem}
\begin{proof}
We focus on constructing the braiding. The remainder of the proof is very similar to that of theorem 3.8 of \cite{DY}. More precisely, let $M$ and $N$ be two local right $B$-modules, and write $$t_{M,N}:M\Box N\rightarrow M\Box^+_B N,\ t_{N,M}:N\Box M\rightarrow N\Box^+_B M$$ for the 2-universal right $B$-balanced 1-morphisms. Then, the composite 1-morphism $t_{N,M}\circ b_{M,N}:M\Box N\rightarrow N\Box^+_B M$ is upgraded to a $B$-balanced right $B$-module 1-morphism via $$\settoheight{\braid}{\includegraphics[width=90mm]{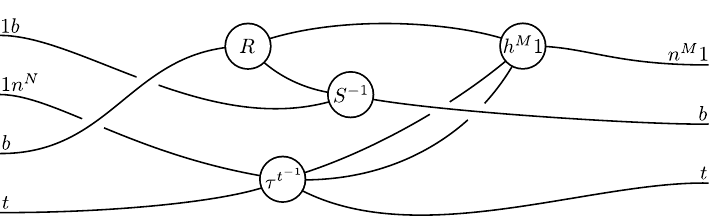}}
\raisebox{0.45\braid}{$\tau^{t\circ b}:=\ $}
\includegraphics[width=90mm]{Pictures/Braiding/tautbholonomy.pdf},$$
$$\settoheight{\braid}{\includegraphics[width=52.5mm]{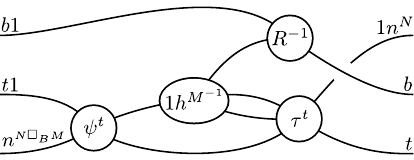}}
\raisebox{0.45\braid}{$\psi^{t\circ b}:=\ $}
\includegraphics[width=52.5mm]{Pictures/Braiding/psitbholonomy.pdf}.$$ Namely, the fact that $\psi^{t\circ b}$ defines a right $B$-module structure on $t_{N,M}\circ b_{M,N}$ follows readily from equation (\ref{eqn:holonomyassociativity}). Further, the proof that $\tau^{t\circ b}$ defines a $B$-balancing follows the argument given in the proof of theorem 3.8 of \cite{DY}, except that we use equation (\ref{eqn:holonomyassociativity}) in lieu of the axioms of their definitions 2.1.2 and 3.2.

Thus, appealing to the the 2-universal property of $t_{M,N}$, the solid arrow diagram below can be filled using a $B$-balanced right $B$-module 2-isomorphism: \begin{equation}\label{eqn:definitionbraiding}\begin{tikzcd}[sep=small]
M\Box N \arrow[dd, "b_{M,N}"'] \arrow[rr, "t_{M,N}"] &  & M\Box^+_B N \arrow[dd, "\widetilde{b}_{M,N}", dotted]  \\
 &  & \\
N\Box M \arrow[rr, "t_{N,M}"'] \arrow[rruu, Rightarrow, "\xi_{M,N}", shorten > = 3ex, shorten < = 4ex]                                   &  & N\Box^+_B M.                                                         
\end{tikzcd}\end{equation} In order to check that $\widetilde{b}$ is compatible with the holonomies, we will appeal to the criterion of lemma \ref{lem:equivalentperspectiveholonomymorphisms}. Firstly, there is a commutative diagram of $B$-$B$-bimodules \begin{equation}\label{eqn:square+}\begin{tikzcd}[sep=small]
Ind^+(M)\Box_B Ind^+(N) \arrow[rr, "b^+"] \arrow[dd, "X^+_{M,N}"'] &  & Ind^+(N)\Box_B Ind^+(M) \arrow[dd, "X^+_{N,M}"] \\
&  &  \\
Ind^+(M\Box_B^+N) \arrow[rr, "Ind^+(\widetilde{b})"']  &  & Ind^+(N\Box_B^+M).      
\end{tikzcd}\end{equation} More precisely, the $B$-$B$-bimodule 1-morphism $b^+$ is defined using a variant of the construction of $\widetilde{b}$, in which the square (\ref{eqn:definitionbraiding}) is viewed in the 2-category of $B$-$B$-bimodules. In order to do this, it is enough to endow $t\circ b$ with a compatible left $B$-module structure via
$$\settoheight{\braid}{\includegraphics[width=52.5mm]{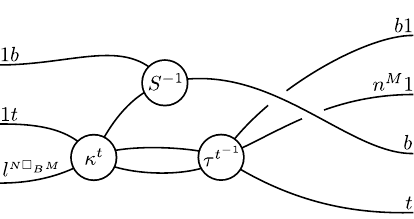}}
\raisebox{0.45\braid}{$\chi^{t\circ b}:=\ $}
\includegraphics[width=52.5mm]{Pictures/Braiding/chitbholonomy.pdf}.$$ It follows easily from the definitions that the square (\ref{eqn:square+}) commutes.

Secondly, also note that a square similar to (\ref{eqn:square+}) with $+$ replaced by $-$ commutes. More precisely, we upgrade the 1-morphism $$t_{N,M}\circ b_{M,N}:Ind^-(M)\Box Ind^-(N)\rightarrow Ind^-(M)\Box_B Ind^-(N)$$ to a right $B$-module $B$-balanced 1-morphism via $$\settoheight{\braid}{\includegraphics[width=90mm]{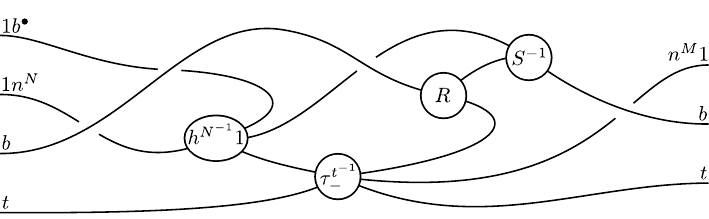}}
\raisebox{0.45\braid}{$\tau_-^{t\circ b}:=\ $}
\includegraphics[width=90mm]{Pictures/Braiding/tautbholonomyvariant.pdf},$$
$$\settoheight{\braid}{\includegraphics[width=52.5mm]{Pictures/Braiding/psitbholonomy.pdf}}
\raisebox{0.45\braid}{$\psi_-^{t\circ b}:=\ $}
\includegraphics[width=52.5mm]{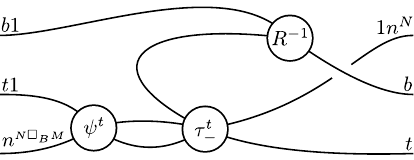}.$$ We claim that this produces the right $B$-module 1-morphism $\widetilde{b}:M\Box_B^-N\rightarrow N\Box_B^-M$. Namely, recall the convention that $M\Box_B^-N = M\Box_B^+N$, which we have adopted in the proof of proposition \ref{prop:2catlocalmodulemonoidal}. In particular, the $B$-balanced structure of $t:M\Box N\rightarrow M\Box_B^-N$ is given by $\tau^t_- = \tau^t \circ (Id\Box (\hbar^N)^{-1})$. The claim then follows by checking that $\widetilde{b}$ satisfies the appropriate universal 2-property. Further, endowing $t_{N,M}\circ b_{M,N}:Ind^-(M)\Box Ind^-(N)\rightarrow Ind^-(N)\Box_B Ind^-(M)$ with a compatible left $B$-module structure $$\settoheight{\braid}{\includegraphics[width=52.5mm]{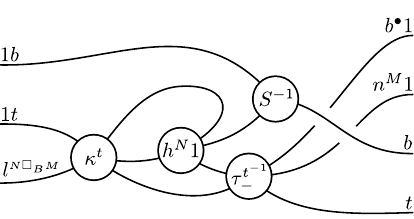}}
\raisebox{0.45\braid}{$\chi_-^{t\circ b}:=\ $}
\includegraphics[width=52.5mm]{Pictures/Braiding/chitbholonomyvariant.pdf}.$$ and appealing to the 2-universal property produces a $B$-$B$-bimodule 1-morphism $b^-$. In summary, we have a commutative square of $B$-$B$-bimodules $$\label{eqn:square-}\begin{tikzcd}[sep=small]
Ind^-(M)\Box_B Ind^+(N) \arrow[rr, "b^-"] \arrow[dd, "X^-_{M,N}"'] &  & Ind^+(N)\Box_B Ind^-(M) \arrow[dd, "X^-_{N,M}"] \\
&  &  \\
Ind^-(M\Box_B^-N) \arrow[rr, "Ind^-(\widetilde{b})"']  &  & Ind^-(N\Box_B^-M).      
\end{tikzcd}$$

Thirdly, we claim that there exists a 2-isomorphism of $B$-$B$-bimodules that makes the diagram below commute $$\begin{tikzcd}[sep=small]
Ind^+(M)\Box_B Ind^+(N) \arrow[dd, "\hbar^M\Box_B \hbar^N"'] \arrow[rr, "b^+"] &  & Ind^+(N)\Box_B Ind^+(M) \arrow[dd, "\hbar^N\Box_B \hbar^M"] \\
 &  & \\
Ind^-(M)\Box_B Ind^-(N) \arrow[rr, "b^-"'] &  & Ind^-(N)\Box_B Ind^-(M).                                  
\end{tikzcd}$$ Thanks to the 2-universal property of the relative tensor product, it is enough to show that there exists a $B$-balanced $B$-$B$-bimodule 2-isomorphism that fits into the next diagram $$\begin{tikzcd}[sep=small]
Ind^+(M)\Box Ind^+(N) \arrow[rr, "b"] \arrow[dd, "\hbar^M\Box \hbar^N"'] &  & Ind^+(N)\Box Ind^+(M) \arrow[rr, "\hbar^N\Box \hbar^M"] &  & Ind^-(N)\Box Ind^-(M) \arrow[dd, "t"] \\
 &  &  &  & \\
Ind^-(M)\Box Ind^-(N) \arrow[rr, "b"']  &  & Ind^-(N)\Box Ind^-(M) \arrow[rr, "t"']  &  & Ind^-(N)\Box_B Ind^-(M).              
\end{tikzcd}$$ We emphasize that $b$ does not carry a $B$-$B$-bimodule structure. Rather, the bottom-left and top-right composite 1-morphisms do. The desired 2-isomorphism is then supplied by the naturality 2-isomorphism of the 2-natural transformation $b$. One checks by tracing through the definitions that this 2-isomorphism is $B$-balanced, and is compatible with the $B$-$B$-bimodule structures.

Putting the above discussion together, we do find that the right $B$-module 1-morphism $\widetilde{b}_{M,N}$ satisfies the criterion of lemma \ref{lem:equivalentperspectiveholonomymorphisms}, so that it is compatible with the holonomies. Appealing to the 2-universal property of the relative tensor product $\Box_B$, it is standard to check that the collection of local right $B$-module 1-morphisms $\widetilde{b}_{M,N}$ for varying $M$ and $N$ assemble into a 2-natural equivalence $\widetilde{b}$.

Finally, analogously to what is done in the proof of theorem 3.8 of \cite{DY}, we can construct invertible modification $\widetilde{R}$ and $\widetilde{S}$ witnessing that $\widetilde{b}$ is appropriately coherent by using the 2-universal property of the relative tensor product. There is no further complication as 2-morphisms of local right $B$-modules are nothing but 2-morphisms of right $B$-modules. This endows $\mathbf{Mod}^{loc}_{\mathfrak{B}}(B)$ with a braiding in the sense of definition 2.3 of \cite{SP} as desired.
\end{proof}

\begin{Remark}\label{rem:laxmonoidalforgetful}
It follows from the proof of the above results that the canonical forgetful 2-functor $$\mathbf{U}:\mathbf{Mod}^{loc}_{\mathfrak{B}}(B)\rightarrow \mathfrak{B}$$ admits a lax braided monoidal structure given on the objects $M$ and $N$ of $\mathbf{Mod}^{loc}_{\mathfrak{B}}(B)$ by $t_{M,N}:M\Box N\rightarrow M\Box^+_B N$, the 1-morphism supplied by the 2-universal property of the relative tensor product over $B$. For later use, let us note that the 2-natural transformation $t$ is strong, and that the relevant coherence modifications as in definition 2.5 of \cite{SP} are all invertible.
\end{Remark}

\begin{Remark}\label{rem:centralfunctor}
The construction of the braiding given in the proof of theorem \ref{thm:2catlocalmodulesbraided} does not use the fact that the $B$-module $N$ is local. In fact, this construction can be upgraded to a braided monoidal 2-functor $$\mathbf{Mod}_{\mathfrak{B}}^{loc}(B)\rightarrow \mathscr{Z}(\mathbf{Mod}_{\mathfrak{B}}(B)).$$
\end{Remark}

\subsection{Local Modules in Braided Fusion 2-Categories}

Let $\mathds{k}$ be an algebraically closed field of characteristic zero. We will now specifically study the properties of the 2-category of local modules in a braided multifusion 2-category $\mathfrak{B}$. We note that the hypotheses on the characteristic of $\mathds{k}$ can be dropped if desired, but we will keep it for simplicity.

\begin{Proposition}\label{prop:finitesemisimple}
Let $B$ be an \'etale algebra in a braided multifusion 2-category $\mathfrak{B}$, then the 2-category $\mathbf{Mod}^{loc}_{\mathfrak{B}}(B)$ of local $B$-modules is finite semisimple.
\end{Proposition}
\begin{proof}
Let $M$ and $N$ be two local $B$-modules in $\mathfrak{B}$. By definition, the 1-category $Hom_B^{loc}(M,N)$ of morphisms of local $B$-modules is a full sub-1-category $Hom_B(M,N)$, which is a finite semisimple 1-category by \cite{D7}. Further, it is easy to check that $Hom_B^{loc}(M,N)$ is closed under direct sums and splittings of idempotents, so that $Hom_B^{loc}(M,N)$ is also a finite semisimple 1-category.

It is clear that the 2-category $\mathbf{Mod}^{loc}_{\mathfrak{B}}(B)$ has direct sums for objects. Next, we show that $\mathbf{Mod}^{loc}_{\mathfrak{B}}(B)$ is closed under the splitting of 2-condensation monads. Take a local $B$-module $M$, together with a 2-condensation monad $(M,e,\xi,\delta)$ in the 2-category $\mathbf{Mod}^{loc}_{\mathfrak{B}}(B)$. By forgetting the holonomy, we obtain a 2-condensation monad on $M$ viewed as a right $B$-module. By proposition 3.3.8 of \cite{D4}, all 2-condensation monads in $\mathbf{Mod}_{\mathfrak{B}}(B)$ split. Thus, there exists a 2-condensation $(M,N,f,g,\phi,\gamma)$ in $\mathbf{Mod}_{\mathfrak{B}}(B)$, together with a splitting $\theta:g \circ f \simeq e$, i.e. a right $B$-module 2-isomorphism satisfying $$ \xi = \theta \cdot (g \circ \phi \circ f) \cdot (\theta^{-1} \circ \theta^{-1}), $$
$$ \delta = (\theta \circ \theta) \cdot (g \circ \gamma \circ f) \cdot \theta^{-1}. $$ We claim that $N$ can be endowed with a holonomoy $h^N$ on $N$ compatible with both $f$ and $g$. This proves that $\mathbf{Mod}^{loc}_{\mathfrak{B}}(B)$ is a Cauchy complete 2-category.

In order to prove the above claim, we let $h^N$ be the 2-isomorphism given by: 
\settoheight{\algebra}{\includegraphics[width=63.75mm]{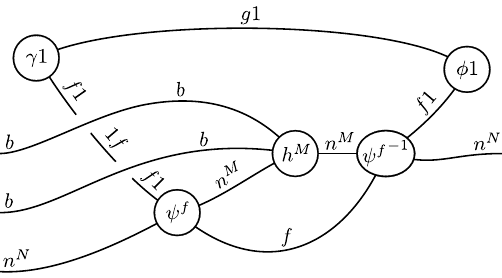}}

\begin{equation}\label{def:holonomycondesation}
    \begin{tabular}{@{}cc@{}}
        \raisebox{0.45\algebra}{$h^N :=$} &
        \includegraphics[width=63.75mm]{Pictures/FiniteSemisimple/definitionholonomycondensation.pdf}.
    \end{tabular}
\end{equation}
        
We now show that $h^N$ does define a holonomy on $N$ using the figures given in appendix \ref{sub:FiniteSemisimpleProof}. The left hand-side of equation (\ref{eqn:holonomyassociativity}) is depicted in figure \ref{fig:cauchycompletenessI1}. To get to figure \ref{fig:cauchycompletenessI2}, we move the coupon labeled $\gamma1$ to the left along the blue arrow. Then, figure \ref{fig:cauchycompletenessI3} is obtained by moving the indicated string to the top along the blue arrows, as well as applying equation (\ref{eqn:moduleassociativity}) to the green coupons. Figure \ref{fig:cauchycompletenessI4} is attained by applying equation (\ref{eqn:holonomyassociativity}) to the blue coupons. So as to obtain figure \ref{fig:cauchycompletenessI5}, we move the coupons labelled $S^{-1}$ and $R^{-1}$ to the left along the blue arrows, followed by moving the coupons labelled $\gamma 11$, $\psi^f1$, and $\phi1$ to the top along the indicated green arrows, then apply equation (\ref{eqn:moduleassociativity}) on the red coupons. Figure \ref{fig:cauchycompletenessI6} is produced by both applying equation (\ref{eqn:module2map}) on the blue coupons and creating a pair of cancelling coupons labeled $\gamma 11$ and $\phi11$ in the green region. We get to figure \ref{fig:cauchycompletenessI7} by moving the freshly created coupon labelled $\gamma11$ to the left along the blue arrow, and inserting a pair of cancelling coupons labelled $\theta11$ and $\theta^{-1}11$ in the green region. Figure \ref{fig:cauchycompletenessI8} is obtained via the use of equation (\ref{eqn:module2map}) on the blue coupons. We then apply equation (\ref{eqn:localmorphism})  on the blue coupons, bringing us to figure \ref{fig:cauchycompletenessI9}. Figure \ref{fig:cauchycompletenessI10} is subsequently obtained by applying equation (\ref{eqn:module2map}) twice. Cancelling the two coupons in blue brings us to figure \ref{fig:cauchycompletenessI11}, which depicts the right hand-side of equation (\ref{eqn:holonomyassociativity}). Further, equation (\ref{eqn:holonomyassociativityprime}) for $h^N$ follows from repeated application of equation (\ref{eqn:modulemapassociativity}) together with equation (\ref{eqn:holonomyassociativityprime}) for $h^M$. Finally, equation (\ref{eqn:holonomyunitality}) for $h^N$ follows from two applications of (\ref{eqn:modulemapunitality}) and one application of equation (\ref{eqn:holonomyunitality}) for $h^M$. Likewise, it is not hard to show that $f$ and $g$ are local $B$-module 1-morphisms.

Using a similar type of argument, one shows that $\mathbf{Mod}^{loc}_{\mathfrak{B}}(B)$ has adjoints for 1-morphisms. Namely, if $f:M\rightarrow N$ is a 1-morphism of local right $B$-modules, then it is not hard to show that the 2-isomorphism

\settoheight{\algebra}{\includegraphics[width=37.5mm]{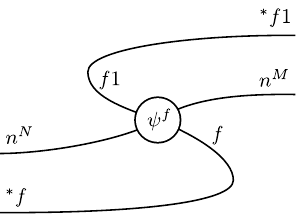}}

\begin{center}
\begin{tabular}{@{}cc@{}}

\raisebox{0.45\algebra}{$\xi^f=$} &
\includegraphics[width=37.5mm]{Pictures/FiniteSemisimple/xif.pdf}
\end{tabular}
\end{center}

\noindent defined in the proof of proposition 2.2.5 of \cite{D7} satisfies the opposite of equation (\ref{eqn:localmorphism}). But, it follows from the same proof that $\xi^f = (\psi^{^*f})^{-1}$, so that the left adjoint of $f$ is a 1-morphism of local right $B$-modules.

Finally, it remains to prove that $\mathbf{Mod}^{loc}_{\mathfrak{B}}(B)$ has finitely many equivalence classes of simple objects. To see this, recall that the forgetful 2-functor $\mathbf{T}:\mathbf{Mod}^{loc}_{\mathfrak{B}}(B)\rightarrow \mathbf{Mod}^+_{\mathfrak{B}}(B)$ is fully faithful on 2-morphisms. In particular, it sends simple objects to simple objects. Thence, it is enough to show that every simple right $B$-module $M$ only admits finitely many holonomies up to equivalence. It follows from lemma \ref{lem:equivalentperspectiveholonomy} that a holonomy on $M$ corresponds exactly to the data of an upgrade of the canonical right $B$-module 1-morphism $Id_M:M\rightarrow M$ to a 1-morphism of $B$-$B$-bimodules $Ind^+(M)\rightarrow Ind^-(M)$. This shows that the set of holonomies on $M$ up to equivalence injects into the set of equivalence classes of invertible $B$-$B$-bimodules 1-morphisms $Ind^+(M)\rightarrow Ind^-(M)$. The later set is finite thanks to theorem 3.1.6 of \cite{D7} and corollary 2.2.3 of \cite{D5}. This finishes the proof.
\end{proof}

\begin{Proposition}\label{prop:duallocalmodules}
Let $\mathfrak{B}$ be a braided multifusion 2-category, and $B$ an \'etale algebra in $\mathfrak{B}$, then $\mathbf{Mod}^{loc}_{\mathfrak{B}}(B)$ has duals.
\end{Proposition}
\begin{proof}
The proof follows by combining lemma \ref{lem:equivalentperspectiveholonomy} with the construction of the dual 2-functor given in appendix \ref{sec:dual}. More precisely, let $M$ be a local right $B$-module in $\mathfrak{B}$. Equivalently, we can think of the holonomy on $M$ as the data necessary to upgrade the canonical right $B$-module 1-morphism $Id_M:M\rightarrow M$ in $\mathfrak{B}$ to a 1-morphism of $B$-$B$-bimodules $\hbar^M:Ind^+(M)\rightarrow Ind^-(M)$. Let us also note that $\hbar^M$ is a 1-isomorphism, and write $(\hbar^M)^{-1}$ for its inverse. Thanks to our assumptions, it follows from \cite{D8} and the proof of theorem 3.6 of \cite{DY} that the monoidal 2-categories $\mathbf{Mod}^+_{\mathfrak{B}}(B)$ and $\mathbf{Bimod}_{\mathfrak{B}}(B)$ have duals. In particular, holonomies on $M^{\sharp}$ correspond to the data of an upgrade of the canonical right $B$-module 1-morphism $Id_{M^{\sharp}}:M^{\sharp}\rightarrow M^{\sharp}$ to a 1-morphism of $B$-$B$-bimodules $Ind^+(M^{\sharp})\rightarrow Ind^-(M^{\sharp})$. By lemma \ref{lem:dual2natural}, up to equivalence, this last $B$-$B$-bimodule 1-morphism is $Ind^+(M)^{\sharp}\rightarrow Ind^-(M)^{\sharp}$. Thus, we can endow $M^{\sharp}$ with the holonomy corresponding to the $B$-$B$-bimodule 1-morphism $\big((\hbar^M)^{-1}\big)^{\sharp}:Ind^+(M^{\sharp})\rightarrow Ind^-(M^{\sharp})$, i.e.\ we have a commuting diagram of $B$-$B$-bimodule 1-morphisms $$\begin{tikzcd}[sep=small]
Ind^+(M^{\sharp}) \arrow[rrr, "\hbar^{M^{\sharp}}"] \arrow[dd, "\rotatebox{90}{$\simeq$}"'] &  &  & Ind^-(M^{\sharp}) \arrow[dd, "\rotatebox{-90}{$\simeq$}"] \\
&  &  & \\
Ind^+(M)^{\sharp} \arrow[rrr,"\big((\hbar^M)^{-1}\big)^{\sharp}"']&  &  & Ind^-(M)^{\sharp}.                                     
\end{tikzcd}$$

It remains to check that ${M^{\sharp}}$ equipped with this holonomy is a right dual for $M$ in $\mathbf{Mod}_{\mathfrak{B}}^{loc}(B)$. It is enough to check that the evaluation and coevaluation 1-morphisms $i_M:B\rightarrow M^{\sharp}\Box^+_B M$ and $e_M:M\Box^+_B M^{\sharp}\rightarrow B$ in $\mathbf{Mod}^+_{\mathfrak{B}}(B)$ are compatible with the holonomies. But, recall from the proof of proposition \ref{prop:2catlocalmodulemonoidal} that the holonomy on the product $M\Box^+_B N$ of two local $B$-module $M$ and $N$ is, up to coherence 1-morphisms, the one corresponding to the $B$-$B$-bimodule 1-morphism $\hbar^M\Box_B \hbar^N$. The claim therefore follows from lemma \ref{lem:equivalentperspectiveholonomymorphisms} by unfolding the definitions and appealing to lemma \ref{lem:dual2natural}.
\end{proof}

\begin{Remark}
More generally, it is not necessary to assume that $\mathfrak{B}$ be multifusion, and that $B$ be separable. The proof of proposition \ref{prop:duallocalmodules} continues to hold provided that $\mathfrak{B}$ has duals, and that the relative tensor product over $B$ of any two $B$-modules exists in $\mathfrak{B}$ and commutes with the monoidal structure.
\end{Remark}

Combining together the two propositions above, we obtain our second main result.

\begin{Theorem}\label{thm:localmodulesmultifusion}
Let $B$ be an \'etale algebra in a braided multifusion 2-category $\mathfrak{B}$. Then, $\mathbf{Mod}_{\mathfrak{B}}^{loc}(B)$ is a braided multifusion 2-category.
\end{Theorem}

\begin{Corollary}
Let $B$ be a connected \'etale algebra in a braided multifusion 2-category $\mathfrak{B}$. Then, $\mathbf{Mod}_{\mathfrak{B}}^{loc}(B)$ is a braided fusion 2-category.
\end{Corollary}

\section{Applications and Examples}

\subsection{Braided Module 1-Categories}

We begin by examining more precisely the notion of local module in the braided fusion 2-category $\mathfrak{B} = \mathbf{2Vect}$. We note that the next results hold more generally without the semisimplicity assumptions, but, for simplicity, we will focus on the semisimple case and work over an algebraically closed field of characteristic zero. Throughout, we work over a braided fusion 1-category $\mathcal{B}$, whose underlying fusion 1-category is assumed to be strict without loss of generality, and with braiding denoted by $\beta$. We will compare the notion of a local $\mathcal{B}$-module in $\mathbf{2Vect}$ with that of a finite semisimple braided $\mathcal{B}$-module 1-category introduced in section 4 of \cite{DN}. We begin by unfolding our definition of a holonomy in the particular case under consideration.

\begin{Definition} \label{def:braidedmodulecat}
A local right $\mathcal{B}$-module 1-category consists of:
\begin{enumerate}
    \item A right $\mathcal{B}$-module 1-category $\mathcal{M}$, with coherence natural isomorphism $\alpha$ given on $M$ in $\mathcal{M}$, and $B$, $C$ in $\mathcal{B}$ by $$\alpha_{M,B,C}:(M\otimes B)\otimes C\rightarrow M\otimes (B\otimes C);$$
    \item A holonomy, that is a natural isomorphism $h$ given on $M$ in $\mathcal{M}$ and $B$ in $\mathcal{B}$ by $$h_{M,B}:M \otimes B \to M \otimes B,$$
\end{enumerate}

satisfying:

\begin{enumerate}
    \item [a.] We have $h_{M,I} = \mathrm{Id}_M $ for all $M$ in $\mathcal{M}$,
    \item [b.] For every $M$ in $\mathcal{M}$ and $B$, $C$ in $\mathcal{B}$, the following diagram commutes:
\end{enumerate}

\begin{equation} \label{eqn:braidedmodulecat1}
   \centering
    \begin{tikzcd}[sep=small]
(M \otimes B) \otimes C \arrow[rrrr, "{h_{M,B} \otimes Id_C}"] \arrow[dd, "{\alpha_{M,B,C}}"'] &  &                                                                     &  & (M \otimes B) \otimes C \arrow[dd, "{h_{M\otimes B,C}}"] \\
 &  &   &  &  \\
M \otimes (B \otimes C) \arrow[dd, "{Id_M\otimes \beta_{B,C}}"']  &  &  &  & (M \otimes B) \otimes C \arrow[dd, "{\alpha_{M,B,C}}"]    \\
 &  &  &  & \\
M \otimes (C \otimes B) \arrow[rrdd, "{Id_M\otimes \beta_{C,B}}"'] &  &  &  & M \otimes (B \otimes C)  \\
 &  &  &  &   \\
 &  & {M \otimes (B \otimes C),} \arrow[rruu, "{h_{M,B\otimes C}}"'] &  & 
\end{tikzcd}
\end{equation}

\begin{enumerate}
    \item [c.] For every $M$ in $\mathcal{M}$ and $B$, $C$ in $\mathcal{B}$, the following diagram commutes:
\end{enumerate}

\begin{equation} \label{eqn:braidedmodulecat2}
    \centering
    \begin{tikzcd}[sep=small]
(M \otimes B) \otimes C \arrow[rrr, "{h_{M \otimes B,C}}"] \arrow[dd, "{\alpha_{M,B,C}}"'] &  &  & (M \otimes B) \otimes C \arrow[dd, "{\alpha_{M,B,C}}"]           \\
 &  &  & \\
M \otimes (B \otimes C) \arrow[dd, "{Id_M\otimes \beta_{B,C}}"']  &  &  & M \otimes (B \otimes C) \arrow[dd, "{Id_M\otimes \beta_{C,B}^{-1}}"] \\
 &  &  &  \\
M \otimes (C \otimes B) \arrow[dd, "{\alpha_{M,C,B}^{-1}}"']                                    &  &  & M \otimes (C \otimes B) \arrow[dd, "{\alpha_{M,C,B}^{-1}}"]      \\
 &  &  & \\
(M \otimes C) \otimes B \arrow[rrr, "{h_{M,C}\otimes Id_B}"']  &  &  & (M \otimes C) \otimes B.                                         
\end{tikzcd}
\end{equation}
\end{Definition}

\begin{Definition} \label{def:braidedmodulefunctor}
A $\mathcal{B}$-module functor $F:\mathcal{M} \to \mathcal{N}$ between two local right $\mathcal{B}$-module 1-categories $\mathcal{M}$ and $\mathcal{N}$ with coherence natural isomorphism $s$ is braided if the diagram below commutes for all $B$ in $\mathcal{B}$ and $M$ in $\mathcal{M}$
\begin{equation} \label{eqn:braidedmodulefunctor}
    \centering
    \begin{tikzcd}[sep=small]
            {B \otimes F(M)}
                \arrow[rr,"h_{B,F(M)}"]
                \arrow[dd,"s_{B,M}"']
            & {}
            & {B \otimes F(M)}
                \arrow[dd,"s_{B,M}"]
            \\ {} & {}
            \\ {F(B \otimes M)}
                \arrow[rr,"F(h_{B,M})"']
            & {}
            & {F(B \otimes M).}
    \end{tikzcd}
\end{equation}
\end{Definition}
 
In section 4 of \cite{DN}, a notion of braiding on a left $\mathcal{B}$-module 1-category was introduced. Further, it is shown therein that the 2-category of finite semisimple braided left $\mathcal{B}$-module 1-categories admits a braided monoidal structure. Let us recall that the underlying monoidal structure is given by the relative tensor product over $\mathcal{B}$ as in example \ref{ex:relativetensor2Vect}. We write $\mathbf{LMod}^{br}(\mathcal{B})$ for this braided monoidal 2-category.

\begin{Proposition} \label{prop:braidedmodulecategoriesequivalenttolocalmodules}
Let $\mathcal{B}$ be a braided fusion 1-category. There is an equivalence of braided monoidal 2-categories $$\mathbf{LMod}^{br}(\mathcal{B}) \simeq \mathbf{Mod}^{loc}(\mathcal{B}).$$
\end{Proposition}
\begin{proof}
Firstly, we have to explain our convention for relating right $\mathcal{B}$-module 1-categories with left $\mathcal{B}$-module 1-categories. For any $B$, $C$ in $\mathcal{B}$ and $M$ in $\mathcal{M}$, we set $B\odot M := M\otimes B$, and use the natural isomorphism \begin{align*}(B\otimes C)\odot M =M\otimes (B\otimes C)&\xrightarrow{Id_M\otimes \beta_{C,B}^{-1}}M\otimes (C\otimes B)\\ & \xrightarrow{\alpha^{-1}_{M,C,B}}(M\otimes C)\otimes B=B\odot (C\odot M)\end{align*} to witness the associativity of the left action. In the notation of the previous section, we are considering the underlying left $\mathcal{B}$-module of $Ind^-(\mathcal{M})$. Secondly, let $h$ be a holonomy on the right $\mathcal{B}$-module 1-category $\mathcal{M}$. Then, $\sigma:=h$ is a $\mathcal{B}$-module braiding on the left $\mathcal{B}$-module 1-category $Ind^-(\mathcal{M})$ in the sense of definition 4.1 of \cite{DN}. In particular, this assignment extends to an equivalence of 2-categories $\mathbf{LMod}^{br}(\mathcal{B})\simeq \mathbf{Mod}^{loc}(\mathcal{B})$.

It remains to show that this equivalence is compatible with the braided monoidal structures. Both of the monoidal structures are enhancements of the relative Deligne tensor product over $\mathcal{B}$. More precisely, let $\mathcal{M}$ and $\mathcal{N}$ be two finite semisimple local right $\mathcal{B}$-module 1-categories. On the one hand, the monoidal structure on $\mathbf{Mod}^{loc}(\mathcal{B})$ is obtained by endowing $\mathcal{M}\boxtimes_{\mathcal{B}}^+\mathcal{N}$ with a holonomy. Given objects $M$ in $\mathcal{M}$ and $N$ in $\mathcal{N}$, recall from example \ref{ex:relativetensor2Vect} that the image of $M\boxtimes N$ under the canonical $\mathcal{B}$-balanced functor $t:\mathcal{M}\boxtimes\mathcal{N}\rightarrow \mathcal{M}\boxtimes^+_{\mathcal{B}}\mathcal{N}$ is given by $t(M\boxtimes N)=(V,\tau_+^{-1})$ with $V$ an object of $\mathcal{M}\boxtimes\mathcal{N}$. Further, the $\mathcal{B}$-balanced structure is witnessed by the natural isomorphism $\tau_+$. It then follows from the proof of proposition \ref{prop:2catlocalmodulemonoidal} that the holonomy on $\mathcal{M}\boxtimes^+_{\mathcal{B}}\mathcal{N}$ is completely characterized by the isomorphism \begin{align*}V\otimes (I\boxtimes B)\xrightarrow{Id\boxtimes h^{\mathcal{N}}}V\otimes (I\boxtimes B)&\xrightarrow{\tau_+}V\otimes (B\boxtimes I)\\ &\xrightarrow{h^{\mathcal{M}}\boxtimes Id}V\otimes (B\boxtimes I)\xrightarrow{\tau_+^{-1}}V\otimes (I\boxtimes B)\end{align*} for every $B$ in $\mathcal{B}$. We emphasize that here we are describing the holonomy on $\mathcal{M}\boxtimes^+_{\mathcal{B}}\mathcal{N}$. In particular, the above expression is the inverse of the one given in the proof of \ref{prop:2catlocalmodulemonoidal} because of lemma \ref{lem:equivalentperspectiveholonomy}.

On the other hand, the monoidal structure on $\mathbf{LMod}^{br}(\mathcal{B})$ is obtained by endowing $Ind^-(\mathcal{M})\boxtimes_{\mathcal{B}}^-Ind^-(\mathcal{N})$ with a $\mathcal{B}$-braiding as in remark 4.13 of \cite{DN}. Given objects $M$ in $Ind^-(\mathcal{M})$ and $N$ in $Ind^-(\mathcal{N})$, it follows from the convention taken in the proof of proposition \ref{prop:2catlocalmodulemonoidal} that the image of $M\boxtimes N$ under the canonical $\mathcal{B}$-balanced functor $t:Ind^-(\mathcal{M})\boxtimes Ind^-(\mathcal{N})\rightarrow Ind^-(\mathcal{M})\boxtimes_{\mathcal{B}}^-Ind^-(\mathcal{N})$ is given by $t(M\boxtimes N)=(V,\tau_+^{-1})$. However, the $\mathcal{B}$-balanced structure is witnessed by the natural isomorphism $\tau_-$. Then, the $\mathcal{B}$-braiding is completely characterized by the isomorphism \begin{align*}V\otimes (B\boxtimes I)\xrightarrow{\tau_-^{-1}}V\otimes (I\boxtimes B)&\xrightarrow{Id\boxtimes\sigma^{\mathcal{N}}}V\otimes (I\boxtimes B)\\ &\xrightarrow{\tau_-}V\otimes (B\boxtimes I)\xrightarrow{\sigma^{\mathcal{M}}\boxtimes Id}V\otimes (B\boxtimes I)\end{align*} for every $M$ in $Ind^-(\mathcal{M})$, $N$ in $Ind^-(\mathcal{N})$, and $B$ in $\mathcal{B}$.

These agree up to the left $\mathcal{B}$-module equivalence $$Ind^-(\mathcal{M}\boxtimes_{\mathcal{B}}^+\mathcal{N}) \rightarrow Ind^-(\mathcal{M})\boxtimes_{\mathcal{B}}^-Ind^-(\mathcal{N})$$ whose underlying functor is characterized by $$V\otimes (I\boxtimes B)\xrightarrow{Id\boxtimes (h^{\mathcal{N}})^{-1}}M\boxtimes V\otimes (I\boxtimes B),$$ and whose left $\mathcal{B}$-module structure is induced by $\tau_+$. More precisely, equation (\ref{eqn:braidedmodulefunctor}) holds with these assignments. This follows from the fact that the following diagram commutes $$\begin{tikzcd}[sep = small]
V\otimes (B\boxtimes I) \arrow[rr, "\tau_+^{-1}"] \arrow[rd, "\tau_-^{-1}"'] & & V\otimes (I\boxtimes B) \arrow[ld, "Id\boxtimes h^{\mathcal{N}}"] \\
& V\otimes (I\boxtimes B), &  
\end{tikzcd}$$ for every $B$ in $\mathcal{B}$. This is a consequence of the convention that $\mathcal{M}\boxtimes_{\mathcal{B}}^+\mathcal{N}=\mathcal{M}\boxtimes_{\mathcal{B}}^-\mathcal{N}$ as right $\mathcal{B}$-module 1-categories, which we have imposed during the proof of proposition \ref{prop:2catlocalmodulemonoidal}. Using this, it is straightforward to check that $\mathbf{LMod}^{br}(\mathcal{B})\simeq \mathbf{Mod}^{loc}(\mathcal{B})$ as monoidal 2-categories by appealing to the 2-universal property of the relative Deligne tensor product.

Finally, the two braidings are defined as follows. On the one hand, the proof of theorem \ref{thm:2catlocalmodulesbraided} constructs a braiding $\mathcal{M}\boxtimes^+_{\mathcal{B}}\mathcal{N}\rightarrow \mathcal{N}\boxtimes^+_{\mathcal{B}}\mathcal{M}$ that is completely characterized by the right $\mathcal{B}$-module isomorphism $$V\otimes(B\boxtimes I)\xrightarrow{\tau_+^{-1}}V\otimes(I\boxtimes B)\xrightarrow{Id\boxtimes h^{\mathcal{M}}} V\otimes(I\boxtimes B)$$ for every $B$ in $\mathcal{B}$. On the other hand, by remark 4.13 of \cite{DN}, the braiding $Ind^-(\mathcal{M})\boxtimes_{\mathcal{B}}^-Ind^-(\mathcal{N})\rightarrow Ind^-(\mathcal{N})\boxtimes_{\mathcal{B}}^-Ind^-(\mathcal{M})$ is characterised by $$V\otimes (B\boxtimes I)\xrightarrow{h^{\mathcal{N}}\boxtimes Id}V\otimes (B\boxtimes I)\xrightarrow{\tau_-^{\dagger}} V\otimes (I\boxtimes B)$$ for every $B$ in $\mathcal{B}$. Here $\tau_-^{\dagger}$ refers to the transpose of $\tau_-$. But, it follows from the 2-universal property of the relative Deligne tensor product that $\tau_-^{\dagger}$ and $\tau_+^{-1}$ are isomorphic $\mathcal{B}$-balancing. Using this along with the 2-universal property of the relative Deligne tensor product, one checks that $\mathbf{LMod}^{br}(\mathcal{B})\simeq \mathbf{Mod}^{loc}(\mathcal{B})$ as braided monoidal 2-categories, which concludes the proof.
\end{proof}

One of the main motivation behind the study of braided $\mathcal{B}$-module 1-categories is that they can be used to model the Drinfeld center of the associated 2-category of $\mathcal{B}$-module 1-categories, as shown in theorem 4.11 of \cite{DN}. In particular, this readily gives the following corollary.

\begin{Corollary}
Let $\mathcal{B}$ be a braided fusion 1-category. Then, we have an equivalence of braided fusion 2-categories $$\mathbf{Mod}^{loc}(\mathcal{B}) \simeq \mathscr{Z}(\mathbf{Mod}(\mathcal{B})).$$
\end{Corollary}

\subsection{Braided Algebras and Local Modules in the 2-Category of Local Modules}

Let us now fix a braided monoidal 2-category $\mathfrak{B}$, and $A$ a braided algebra in $\mathfrak{B}$ for which the relative tensor product of any two modules exists and commutes with the monoidal product. We have seen in theorem \ref{thm:2catlocalmodulesbraided} that the 2-category $\mathbf{Mod}_{\mathfrak{B}}^{loc}(A)$ is braided monoidal. It is therefore natural to ask what are the braided algebras in this 2-category.
    
\begin{Lemma}
The data of a braided algebra in $\mathbf{Mod}^{loc}_\mathfrak{B}(A)$ corresponds exactly to the data of a braided algebra $B$ in $\mathfrak{B}$ equipped with a 1-homomorphism of braided algebras $f:A \rightarrow B$ in $\mathfrak{B}$.
\end{Lemma}
\begin{proof}
It follows immediately from remark \ref{rem:laxmonoidalforgetful} that the data of a braided algebra in $\mathbf{Mod}^{loc}_\mathfrak{B}(A)$ gives a braided algebra in $\mathfrak{B}$ equipped with a 1-homomorphism of braided algebras from $A$. Conversely, given any braided algebra $B$ in $\mathfrak{B}$ equipped with a 1-homomorphism of braided algebras $f:A \rightarrow B$. Then, the braiding on $B$ yields a canonical holonomy on $B$ viewed as a right $A$-module. Further, it follows from the 2-universal property of the relative tensor product over $A$ that the multiplication 1-morphism $m^B:B\Box B\rightarrow B$ factors as $$B\Box B\xrightarrow{t_{B,B}}B\Box^+_A B\xrightarrow{\hat{m}}B.$$ It is not difficult to check that $\hat{m}$ is compatible with the canonical holonomy on $B$,  the remaining data as well as the coherence conditions for a braided algebra in $\mathbf{Mod}^{loc}_\mathfrak{B}(A)$ follow from those of $B$ via the 2-universal property of $\Box^+_A$. This gives the desired result.
\end{proof}

Let us now fix a 1-homomorphism of braided algebras $f:A \rightarrow B$ in $\mathfrak{B}$, and assume that the relative tensor products over both $A$ and $B$ exist in $\mathfrak{B}$ and commute with the monoidal structure. Using the underlying 1-homomorphism of algebras $f:A\rightarrow B$, we can view $B$ as an algebra in $\mathbf{Bimod}_{\mathfrak{B}}(A)$. Our hypotheses guarantee that the relative tensor product over $B$ exists in $\mathbf{Bimod}_{\mathfrak{B}}(A)$, and is in fact given by the relative tensor product over $B$ in $\mathfrak{B}$ (see corollary 3.2.12 of \cite{D8}). More precisely, the canonical lax monoidal 2-functor $$\mathbf{Bimod}_{\mathbf{Bimod}_{\mathfrak{B}}(A)}(B)\rightarrow \mathbf{Bimod}_{\mathfrak{B}}(B)$$ is in fact strongly monoidal as well as an equivalence. In particular, we find that the canonical lax monoidal 2-functor $$\mathbf{Mod}^+_{\mathbf{Mod}^+_{\mathfrak{B}}(A)}(B)\rightarrow \mathbf{Mod}^+_{\mathfrak{B}}(B)$$ is a (strongly) monoidal equivalence. Then, it is also natural to ask how the 2-category of local $B$-modules in $\mathfrak{B}$ and in $\mathbf{Mod}^{loc}_\mathfrak{B}(A)$ compare. We will see some applications of this result in the next section.

\begin{Proposition}\label{prop:iteratedlocalmodules}
There is an equivalence of braided monoidal 2-categories $$\mathbf{Mod}^{loc}_{\mathbf{Mod}^{loc}_\mathfrak{B}(A)}(B) \simeq \mathbf{Mod}^{loc}_\mathfrak{B}(B).$$
\end{Proposition}
\begin{proof}
Note that there is a lax braided monoidal forgetful 2-functor $$\mathbf{V}:\mathbf{Mod}^{loc}_{\mathbf{Mod}^{loc}_\mathfrak{B}(A)}(B) \rightarrow \mathbf{Mod}^{loc}_\mathfrak{B}(B).$$ One way to see this is to recall from remark \ref{rem:laxmonoidalforgetful} that the forgetful 2-functor $\mathbf{U}:\mathbf{Mod}^{loc}_\mathfrak{B}(A)\rightarrow\mathfrak{B}$ has a lax braided monoidal structure. More precisely, in the notation of definition 2.5 of \cite{SP}, we have strong 2-natural transformations $\chi^{\mathbf{U}}$ and $\iota^{\mathbf{U}}$, which are not necessarily equivalences, and invertible modifications $\omega^{\mathbf{U}}$, $\gamma^{\mathbf{U}}$, $\delta^{\mathbf{U}}$, and $u^{\mathbf{U}}$. Using $f:A\rightarrow B$, we may view $B$ as a braided algebra in $\mathbf{Mod}^{loc}_\mathfrak{B}(A)$, and we manifestly have $\mathbf{U}(B)=B$. Then, we can take local modules over $B$ in both $\mathbf{Mod}^{loc}_\mathfrak{B}(A)$ and $\mathfrak{B}$. But, taking local modules is functorial, so that the lax braided monoidal 2-functor $\mathbf{U}$ induces the lax braided monoidal 2-functor $\mathbf{V}$ upon taking local modules over $B$.

We now show that $\mathbf{V}$ induces an equivalence of 2-categories. Namely, by inspecting the definitions, we find that $\mathbf{V}$ is fully faithful on 2-morphisms. Then, it is not difficult to check that it is essentially surjective on 1-morphisms. The fact that $\mathbf{V}$ is essentially surjective on objects follows from a slight generalization of the argument used in the proof of the lemma above. We leave the details to the keen reader.

It remains to check that $\mathbf{V}$ is compatible with the braided monoidal structures. As the 2-functor $\mathbf{V}$ inherits a canonical lax braided monoidal structure from $\mathbf{U}$, it is enough to show that this lax structure is actually strong. In fact, as all the modifications witnessing the coherence of $\mathbf{U}$ are invertible, this is also true of the modifcations witnessing the coherence of $\mathbf{V}$. Thus, it is enough to show that the (necessarily strong) 2-natural transformations $\chi^{\mathbf{V}}$ and $\iota^{\mathbf{V}}$ are equivalences. To see this, note that there is a commutative square of lax monoidal 2-functors $$\begin{tikzcd}[sep=small]
\mathbf{Mod}^{loc}_{\mathbf{Mod}^{loc}_\mathfrak{B}(A)}(B) \arrow[dd] \arrow[rr, "\mathbf{V}"] &  & \mathbf{Mod}^{loc}_\mathfrak{B}(B) \arrow[dd] \\
&  &\\
\mathbf{Mod}^{+}_{\mathbf{Mod}^{+}_\mathfrak{B}(A)}(B) \arrow[rr, "\simeq"]                    &  & \mathbf{Mod}^{+}_\mathfrak{B}(B).             
\end{tikzcd}$$ But, it follows from the discussion preceding the present lemma that the bottom horizontal arrow is strongly monoidal and an equivalence. Further, we have seen in remark \ref{rem:monoidalforgetful} that the vertical arrows are strongly monoidal, conservative on 1-morphisms, and fully faithful on 2-morphisms. This proves that $\chi^{\mathbf{V}}$ and $\iota^{\mathbf{V}}$ are 2-natural equivalences, so that the lax braided monoidal structure of $\mathbf{V}$ is actually strong as desired.
\end{proof}
    
\subsection{Lagrangian Algebras}

We work over an algebraically closed field of characteristic zero, and fix $\mathfrak{B}$ a braided fusion 2-category 

\begin{Definition}
A Lagrangian algebra in $\mathfrak{B}$ is a connected étale algebra $B$ in $\mathfrak{B}$ such that $\mathbf{Mod}^{loc}_\mathfrak{B}(B) \simeq \mathbf{2Vect}$.
\end{Definition}

\begin{Example}
A Lagrangian algebra in $\mathbf{2Vect}$ is a non-degenerate braided fusion 1-category, as can be seen from proposition \ref{prop:braidedmodulecategoriesequivalenttolocalmodules} and proposition 4.17 of \cite{DN}.
\end{Example}
    
\begin{Remark}
It follows from the previous example that the property of being Lagrangian for a connected étale algebra is a categorical non-degeneracy condition.
\end{Remark}

Recall that connected étale algebras in $\mathbf{2Vect}$ are braided fusion 1-categories. Further, given a braided fusion 1-category $\mathcal{B}$, we have seen in corollary \ref{cor:LagrangianalgebrainZModB} that $\mathbf{Mod}^{loc}(\mathcal{B}) \simeq \mathscr{Z}(\mathbf{Mod}(\mathcal{B}))$. The next corollary therefore follows by proposition \ref{prop:iteratedlocalmodules}.
    
\begin{Corollary}\label{cor:LagrangianalgebrainZModB}
Let $\mathcal{B}$ be any braided fusion 1-category. Lagrangian algebras in $\mathscr{Z}(\mathbf{Mod}(\mathcal{B}))$ correspond exactly to non-degenerate braided fusion 1-categories $\mathcal{C}$ equipped with a braided monoidal functor $\mathcal{B} \rightarrow \mathcal{C}$.
\end{Corollary}

Another notion of Lagrangian algebra was introduced in section 2.3 of \cite{JFR}. We recall their definition below, and use the name alter-Lagrangian algebras to refer to such objects.

\begin{Definition}
An alter-Lagrangian algebra in a braided fusion 2-category $\mathfrak{B}$ is an étale algebra $B$ in $\mathfrak{B}$ satisfying:

\begin{enumerate}
    \item[a.] It is strongly connected, i.e.\ the 1-morphism $i:I\rightarrow B$ is the inclusion of a simple summand.
    \item[b.] Its M\"uger center is trivial, i.e.\ the fusion 1-category $Hom_B^{loc}(B,B)$ of local right $B$-module 1-morphisms $B\rightarrow B$, is trivial.
\end{enumerate}
\end{Definition}

\begin{Example}
As explained in remark 2.27 of \cite{JFR}, the condition of being strongly connected in the definition of an alter-Lagrangian algebra is simply too strong. This is why we have only insisted that Lagrangian algebras are connected. For instance, let $G$ be a finite group, and consider the braided fusion 2-category $\mathfrak{B}=\mathscr{Z}(\mathbf{Mod}(\mathbf{Rep}(G)))$. We have seen that Lagrangian algebras in $\mathfrak{B}$ correspond exactly to braided fusion 1-categories $\mathcal{C}$ equipped with a braided functor $F:\mathbf{Rep}(G)\rightarrow\mathcal{C}$. The corresponding Lagrangian algebra is strongly connected if and only if $F$ is fully faithful, in which case $F$ is an inclusion. In particular, $\mathbf{Vect}$ equipped with the forgetful functor $\mathbf{Rep}(G)\rightarrow\mathbf{Vect}$ defines a Lagrangian algebra in $\mathscr{Z}(\mathbf{Mod}(\mathbf{Rep}(G)))$ that is not strongly connected, and therefore not an alter-Lagrangian algebra.
\end{Example}

\begin{Proposition}
Let $\mathcal{B}$ be a braided fusion 1-category. Then, every alter-Lagrangian algebra in $\mathscr{Z}(\mathbf{Mod}(\mathcal{B}))$ is Lagrangian.
\end{Proposition}
\begin{proof}
Let $\mathcal{C}$ be a braided fusion 1-category equipped with a fully faithful braided monoidal functor $F:\mathcal{B}\rightarrow \mathcal{C}$. Then, it follows from proposition 2.28 of \cite{JFR} that its M\"uger $\mathcal{Z}_{(2)}(\mathcal{C})$ is trivial, i.e.\ we have $\mathcal{Z}_{(2)}(\mathcal{C})\simeq\mathbf{Vect}$, or, equivalently, $\mathcal{C}$ is non-degenerate. The result then follows from proposition \ref{prop:braidedmodulecategoriesequivalenttolocalmodules}.
\end{proof}

\begin{Remark}
In fact, as was already noted in remark 2.27 of \cite{JFR}, the above proposition should hold for any non-degenerate braided fusion 2-category $\mathfrak{B}$, that is, braided fusion 2-category with trivial sylleptic center in the sense of \cite{Cr}. However, we point out that this is not true for an arbitrary braided fusion 2-category, as can be seen from the example below.
\end{Remark}

\begin{Example}
Let $\mathfrak{B}:=\mathbf{2Vect}_A$ the braided fusion 2-category of $A$-graded 2-vector spaces for some finite abelian group $A$. Then, alter-Lagrangian algebras in $\mathbf{2Vect}_A$ are $A$-graded braided fusion 1-categories whose $0$-graded part is a non-degenerate braided fusion 1-category. In particular, we can view any non-degenerate braided fusion 1-category $\mathcal{B}$ as a connected étale algebra $B$ in $\mathbf{2Vect}_A$, and we have $\mathbf{Mod}^{loc}_{\mathbf{2Vect}_A}(B)\simeq \mathbf{2Vect}_A$. In particular, $B$ is not a Lagrangian algebra in $\mathbf{2Vect}_A$.
\end{Example}

\subsubsection{Lagrangian Algebras in Physics}

We would now like to discuss the significance of Lagrangian algebras in Physics. Using some bootstrap analysis \cite{KWZ,JF}, (3+1)-dimensional topological phases of matter are believed to be characterized by non-degenerate braided fusion 2-categories, which can be thought of as a collection of topological excitations of all codimensions \cite{KWZ} or low-energy topological sectors of observables \cite{JF} in the given quantum many-body systems. It is then natural to ask what mathematical structure corresponds to the (2+1)-dimensional topological boundary conditions of such a (3+1)d topological phase in the bulk.\footnote{We emphasize that we require that a topological boundary satisfy a boundary-bulk correspondence \cite{KW}, i.e.\ that the Drinfeld center of the boundary is equivalent to the bulk. In other word, the condensed bulk phase consisting of deconfined particles is trivial.} Let us write $\mathfrak{Z}$ for the non-degenerate braided fusion 2-category $\mathscr{Z}(\mathfrak{C})$ for some fusion 2-category $\mathfrak{C}$. Generalizing ideas from \cite{K14}, we find that the (2+1)d topological boundary conditions for the (3+1)d phase corresponding to $\mathfrak{Z}$ are given by Lagrangian algebras in $\mathfrak{Z}$. In fact, given a Lagrangian algebra $A$ in $\mathfrak{Z}$, the fusion 2-category $\mathbf{Mod}_{\mathfrak{Z}}(A)$ describes the collection of topological excitations on the boundary, and the canonical braided 2-functor $\mathfrak{Z}\rightarrow \mathscr{Z}(\mathbf{Mod}_{\mathfrak{Z}}(A))$ encodes the interaction between the bulk and the boundary.

To examine the validity of this assertion, it is instructive to consider the case when the bulk is the trivial (3+1)d topological phase, mathematically described by $\mathbf{2Vect}$. By arguments from \cite{KW}, the bulk phase controls all gravitational anomalies of its topological boundary conditions. Thus (2+1)d topological boundary conditions for the trivial (3+1)d topological phase are nothing else but (2+1)d topological phases, which are classified by non-degenerate braided fusion 1-categories as established in \cite{Kit03} for point-like excitations and in \cite{KWZ,KLWZZ,JF} including all string-like excitations. Meanwhile, by corollary \ref{cor:LagrangianalgebrainZModB}, we see that Lagrangian algebras in $\mathbf{2Vect}$ are exactly non-degenerate braided fusion 1-categories. Here mathematics matches perfectly with the physical intuition.

In the case of the (3+1)d toric code model \cite{ZLZHKT}, explicit computations were carried out using a microscopic realization on a $3$d cubical lattice (see also \cite{Luo22}). Microscopically, they discovered three Lagrangian algebras $A_e,A_1,A_2$ in $\mathbf{TC} := \mathscr{Z}(\mathbf{2Vect}_{\mathbb{Z}/2})$, corresponding to a rough boundary condition $\mathbf{Mod}_{\mathbf{TC}}(A_e)$, a smooth boundary condition $\mathbf{Mod}_{\mathbf{TC}}(A_1)$ and a twisted smooth boundary condition $\mathbf{Mod}_{\mathbf{TC}}(A_2)$. We note that $\mathbf{Mod}_{\mathbf{TC}}(A_e)\simeq\mathbf{2Vect}_{\mathbb{Z}/2}$, and $\mathbf{Mod}_{\mathbf{TC}}(A_2)\simeq\mathbf{Mod}_{\mathbf{TC}}(A_1)\simeq\mathbf{2Rep}(\mathbb{Z}/2)$. However, the boundaries provided by $A_1$ and $A_2$ are distinct as the braided 2-functors $\mathbf{TC}\rightarrow \mathscr{Z}(\mathbf{2Rep}(\mathbb{Z}/2))$ are distinct. From these three elementary boundary conditions, one can construct infinitely many others by stacking with an anomaly-free (2+1)d topological phase. In mathematical language, this means that, given a Lagrangian algebra $L$ in $\mathfrak{Z}$, and a Lagrangian algebra $A$ in $\mathbf{2Vect}$, that is a non-degenerate braided fusion 1-category, we obtain a new Lagrangian algebra $L \boxtimes A$ in $\mathfrak{Z}\boxtimes \mathbf{2Vect}\simeq\mathfrak{Z}$.

\begin{figure}[!htbp]
    \centering
    \begin{tikzpicture}[scale=0.8]
        \fill[gray!30] (-4,0) rectangle (0,4) ;
        \fill[gray!20] (0,0) rectangle (5,4) ;
        \fill[gray!10] (5,0) rectangle (10,4) ;
        \draw[dotted] (0,4)--(0,0) node[midway] {$Id$} ;
        \draw[dotted] (5,4)--(5,0) node[midway] {\text{gauging}} ;
        \draw[very thick] (-4,0)--(10,0) ;
        \draw[fill=white] (-0.1,-0.1) rectangle (0.1,0.1) ;
        \draw[fill=white] (4.9,-0.1) rectangle (5.1,0.1) ;
        \node at (0,-0.5) {$\mathbf{2Vect}$} ;
        \node at (5,-0.5) {\text{gauging}} ;
        \node at (-2,2.5) {$\mathscr{Z}(\mathbf{2Vect}_G)$} ;
        \node at (2.7,2.5) {$\mathscr{Z}(\mathbf{2Rep}(G))$} ;
        \node at (7.5,2.5) {$\mathscr{Z}(\mathbf{Mod}(\mathbf{Rep}(G)))$} ;
        \node at (-2,-0.5) {$\mathbf{2Vect}_G$} ;
        \node at (2.7,-0.5) {$\mathbf{2Rep}(G)$} ;
        \node at (7.5,-0.5) {$\mathbf{Mod}(\mathbf{Rep}(G))$} ;
    \end{tikzpicture}
    \caption{Braided equivalences of Drinfeld centers}
    \label{fig:EquivalentThreeCenters}
    \end{figure}
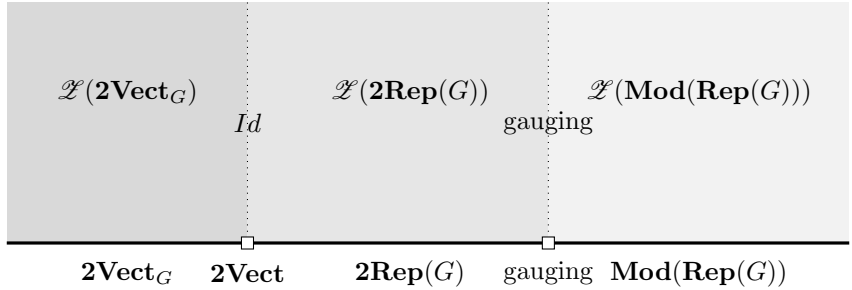

At this point, it is natural to ask whether the above construction exhausts all the possible Lagrangian algebras in the (3+1)d toric code model. As a first step towards answering this question, let us fix a finite group $G$, and consider the braided equivalences between Drinfeld centers depicted in figure \ref{fig:EquivalentThreeCenters}. The two braided fusion 2-categories $\mathscr{Z}(\mathbf{2Vect}_G)$ and $\mathscr{Z}(\mathbf{2Rep}(G))$ are both modelled by the 2-category of finite semisimple $G$-crossed 1-categories \cite{KTZ}, and this is witnessed by the Morita equivalence between fusion 2-categories $\mathbf{2Vect}_G$ and $\mathbf{2Rep}(G)$ \cite{D7, D8}. The equivalence between the braided fusion 2-categories $\mathscr{Z}(\mathbf{2Rep}(G))$ and $\mathscr{Z}(\mathbf{Mod}(\mathbf{Rep}(G)))$ is induced by an equivalence of symmetric fusion 2-categories between $\mathbf{2Rep}(G)$ and $\mathbf{Mod}(\mathbf{Rep}(G))$. More precisely, this equivalence is implemented by equivariantization for finite semisimple 1-categories with a $G$-action \cite{DGNO}. This procedure can be reversed via de-equivariantization, so no information is lost in these processes. In physical language, equivariantization corresponds to gauging a $G$-symmetry on a system, thereby obtaining an equivalent system, which is now equipped with the dual symmetry $\mathbf{Rep}(G)$. If $G$ is abelian, the dual symmetry is given by $\mathbf{Rep}(G)\simeq \mathbf{Vect}_{\widehat{G}}$, which is invertible, but this is not the case in general. %Conversely, de-equivariantization is a form of generalized ``anyon condensation''.

Furthermore, the equivalence $\mathscr{Z}(\mathbf{2Rep}(G))\simeq \mathscr{Z}(\mathbf{Mod}(\mathbf{Rep}(G)))$ is realized by equivariantization for finite semisimple $G$-crossed 1-categories, which produces a finite semisimple braided $\mathbf{Mod}(\mathbf{Rep}(G))$-module 1-category. This process may also be reversed via de-equivariantization. In terms of étale algebras, we obtain an equivalence between $G$-crossed braided multifusion 1-categories, and braided multifusion 1-category equipped with a braided functor from $\mathbf{Rep}(G)$. This correspondence is well-known in the theory of braided fusion 1-categories \cite{Ki, Mu1}.

Then, thanks to corollary \ref{cor:LagrangianalgebrainZModB}, Lagrangian algebras in $\mathscr{Z}(\mathbf{Mod}(\mathbf{Rep}(G)))$ are non-degenerate braided fusion 1-categories equipped with a braided functor from $\mathbf{Rep}(G)$. Fixing such a 1-category $\mathcal{B}$, notice that the braided functor $\mathbf{Rep}(G) \to \mathcal{B}$ must factorize as $$\mathbf{Rep}(G) \twoheadrightarrow \mathbf{Rep}(H) \hookrightarrow \mathcal{B}$$ for some subgroup $H\subseteq G$ (determined up to conjugation). More precisely, the first braided functor is dominant or surjective, whereas the second is fully faithful or injective. In particular, we can view $\mathcal{B}$ as a non-degenerate extension of $\mathbf{Rep}(H)$.

For example, as explained in remark 4.1 of \cite{ZLZHKT}, in the (3+1)d toric code model $\mathbf{TC}$, i.e.\ with $G=\mathbb{Z}/2$, we find that the three Lagrangian algebras $A_e, A_1, A_2$ can be described as follows: \begin{itemize}
    \item The Lagrangian algebra $A_e$ corresponds to the forgetful functor $\mathbf{Rep}(\mathbb{Z}_2) \twoheadrightarrow \mathbf{Vect}$.

    \item The Lagrangian algebra $A_1$ corresponds to the minimal non-degenerate extension $\mathbf{Rep}(\mathbb{Z}_2) \hookrightarrow \mathcal{Z}(\mathbf{Vect}_{\mathbb{Z}/2})$.

    \item The Lagrangian algebra $A_2$ corresponds to the minimal non-degenerate extension $\mathbf{Rep}(\mathbb{Z}_2) \hookrightarrow \mathcal{Z}(\mathbf{Vect}^{\omega}_{\mathbb{Z}/2})$, where $\omega$ is a cocycle representing the non-trivial element in $\mathrm{H}^3(\mathbb{Z}/2;\mathbb{C}^\times) \simeq \mathbb{Z}/2$.
\end{itemize}

More generally, given a minimal non-degenerate extension $\mathbf{Rep}(G) \hookrightarrow \mathcal{M}$ viewed as a Lagrangian algebra in $\mathscr{Z}(\mathbf{Mod}(\mathbf{Rep}(G)))$, and a non-degenerate braided fusion 1-category $\mathcal{A}$ viewed as a Lagrangian algebra in $\mathbf{2Vect}$, we can consider the Lagrangian algebra in $\mathscr{Z}(\mathbf{Mod}(\mathbf{Rep}(G)))$ given by $\mathbf{Rep}(G) \hookrightarrow\mathcal{M} \boxtimes \mathcal{A}$. It is clear that not all non-degenerate extensions of $\mathbf{Rep}(G)$ are of the this form.

In section 4 of \cite{ZLZHKT}, it is proposed that all possible (2+1)d topological boundary conditions for the (3+1)d toric code model are given by first stacking the rough or smooth boundary with an anomaly-free (2+1)d topological order, and then introduce a twist or coupling between them, thereby obtaining a new boundary of $\mathbf{TC}$ from its condensation. Meanwhile, using the theory of Lagrangian algebras we have developed, these boundary conditions corresponds to de-equivariantization of either non-degenerate extensions $\mathbf{Rep}(\mathbb{Z}_2) \hookrightarrow \mathcal{M}$, or forgetful functor $\mathbf{Rep}(\mathbb{Z}_2) \twoheadrightarrow \mathbf{Vect}$ stacked with a non-degenerate braided fusion 1-category $\mathcal{B}$. In the second case, the boundary condition can always be obtained by stacking the Lagrangian algebra $A_e$ in $\mathbf{TC}$ with a non-degenerate braided fusion 1-category. However, in the first case, this is not true. We believe it is an interesting question for physicists to explicitly realize the correspondence between Lagrangian algebras in $\mathscr{Z}(\mathbf{Mod}(\mathbf{Rep}(G)))$ and topological boundary conditions of $\mathbf{TC}$ in the microscopic lattice model. % Old

\appendix

\section{Appendices}

\subsection{Monoidality of the Dual 2-Functor}\label{sec:dual}

Let $\mathfrak{C}$ be a monoidal 2-category, which we will assume without loss of generality is strict cubical. Let us suppose that $\mathfrak{C}$ has right duals, that is, for every object $C$ of $\mathfrak{C}$, there exists an object $C^{\sharp}$ of $\mathfrak{C}$ together with 1-morphisms $i_C:I\rightarrow C^{\sharp}\Box C$, and $e_C:C\Box C^{\sharp}\rightarrow I$ that satisfy the snake equations up to 2-isomorphisms. In this case, it follows from corollary 2.8 of \cite{Pstr} that every objects has a coherent right duals, that is, there exists 2-isomorphisms $$E_C:(e_C\Box C)\circ (C\Box i_C)\Rightarrow Id_C,$$ $$F_C:(C^{\sharp}\Box e_C)\circ (i_C\Box C^{\sharp})\Rightarrow Id_{C^{\sharp}}$$ such that 

\newlength{\adj}
\settoheight{\adj}{\includegraphics[width=37.5mm]{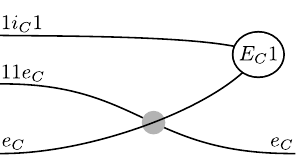}}

\begin{equation}\label{eqn:swallowtail1}
\begin{tabular}{@{}ccc@{}}
\includegraphics[width=37.5mm]{Pictures/Appendix/coherence1.pdf} & \raisebox{0.45\adj}{$=$} &

\includegraphics[width=30mm]{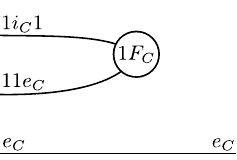}.
\end{tabular}
\end{equation}

\begin{equation}\label{eqn:swallowtail2}
\begin{tabular}{@{}ccc@{}}
\includegraphics[width=37.5mm]{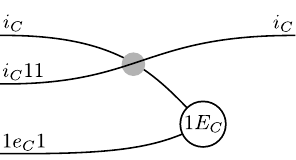} & \raisebox{0.45\adj}{$=$} &

\includegraphics[width=30mm]{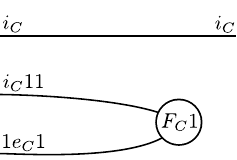}.
\end{tabular}
\end{equation}

\noindent It then follows from appendix A.2 of \cite{D2} that there is a 2-functor $(-)^{\sharp}:\mathfrak{C}\rightarrow \mathfrak{C}^{1op}$ that sends an object to its right dual. Here, we use $\mathfrak{C}^{1op}$ to denote the 2-category obtained from $\mathfrak{C}$ by reversing the direction of the 1-morphisms. We recall the construction of the dual 2-functor in more detail. For any object $C$ of $\mathfrak{C}$, we fix a coherent right dual $(C, C^{\sharp},i_C,e_C,E_C,F_C)$. Then, the 2-functor $(-)^{\sharp}$ sends the object $C$ of $\mathfrak{C}$ to $C^{\sharp}$. Next, a 1-morphism $f:C\rightarrow D$ is sent via $(-)^{\sharp}$ to the 1-morphism in $\mathfrak{C}$ given by $$f^{\sharp}:=(C^{\sharp}\Box e_D)\circ (C^{\sharp}\Box f\Box D^{\sharp})\circ (i_C\Box D^{\sharp}):D^{\sharp}\rightarrow C^{\sharp},$$ and viewed as a 1-morphism in $\mathfrak{C}^{1op}$. Further, a 2-morphism $\alpha:f\Rightarrow g:C\rightrightarrows D$ is sent by $(-)^{\sharp}$ to $$(C^{\sharp}\Box e_D)\circ (C^{\sharp}\Box \alpha\Box D^{\sharp})\circ (i_C\Box D^{\sharp}):f^{\sharp}\Rightarrow g^{\sharp}.$$ Finally, the coherence 2-isomorphisms witnessing unitality for the 2-functor $(-)^{\sharp}$ is given on the object $C$ in $\mathfrak{C}$ by $$\phi^{\sharp}_C:=D_C,$$ and the coherence 2-isomorphisms witnessing that the 2-functor $(-)^{\sharp}$ is compatible with composition of 1-morphisms is given on $f:B\rightarrow C$ and $g:C\rightarrow D$ by

\settoheight{\adj}{\includegraphics[width=60mm]{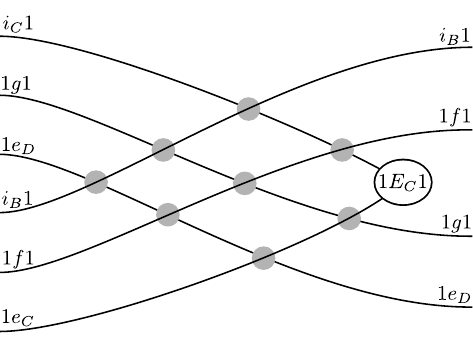}}

$$\begin{tabular}{@{}cc@{}}
\raisebox{0.45\adj}{$\phi^{\sharp}_{f,g}:=$} & \includegraphics[width=60mm]{Pictures/Appendix/phifg.pdf}.
\end{tabular}$$

The result below positively answers a question raised in section 1.2 of \cite{D2}. We use $\mathfrak{C}^{mop,1op}$ to denote the monoidal 2-category obtained from $\mathfrak{C}$ by both taking the opposite monoidal product and reversing the direction of 1-morphisms.

\begin{Proposition}\label{prop:dualfunctorial}
Let $\mathfrak{C}$ be a monoidal 2-category that has right duals. Then, the 2-functor $(-)^{\sharp}:\mathfrak{C}\rightarrow \mathfrak{C}^{mop, 1op}$ sending an object to its right dual admits a canonical monoidal structure.
\end{Proposition}
\begin{proof}
We note that, as $\mathfrak{C}$ is strict cubical, then $\mathfrak{C}^{mop, 1op}$ is strict cubical. We follow the notations of \cite{SP} for the coherence data for $(-)^{\sharp}$. We begin by constructing a 1-equivalence $\iota$ witnessing that $(-)^{\sharp}$ is compatible with the monoidal units. Without loss of generality, we will assume that the chosen coherent right dual to $I$ is $I$ with the identity 1-morphisms and 2-morphisms. In particular, we can take $\iota$ to be the identity 2-natural transformation. We also construct the 2-natural equivalence $\chi$ witnessing that $(-)^{\sharp}$ is compatible with the monoidal structures. More precisely, given objects $C$ and $D$ in $\mathfrak{C}$, we define a 1-morphism $(C\Box D)^{\sharp}\rightarrow D^{\sharp}\Box C^{\sharp}$ in $\mathfrak{C}$ by $$\chi_{C,D}:=(D^{\sharp}\Box C^{\sharp}\Box e_{C\Box D})\circ (D^{\sharp}\Box i_{C} \Box D\Box (C\Box D)^{\sharp})\circ (i_{D}\Box (C\Box D)^{\sharp}),$$ which, when viewed as a 1-morphism in $\mathfrak{C}^{mop,1op}$, provides the underlying 1-morphisms of the requisite coherence 2-natural transformation. The 2-naturality of $\chi$ is given on the 1-morphisms $f:A\rightarrow B$ and $g:C\rightarrow D$ in $\mathfrak{C}$ by

\settoheight{\adj}{\includegraphics[width=90mm]{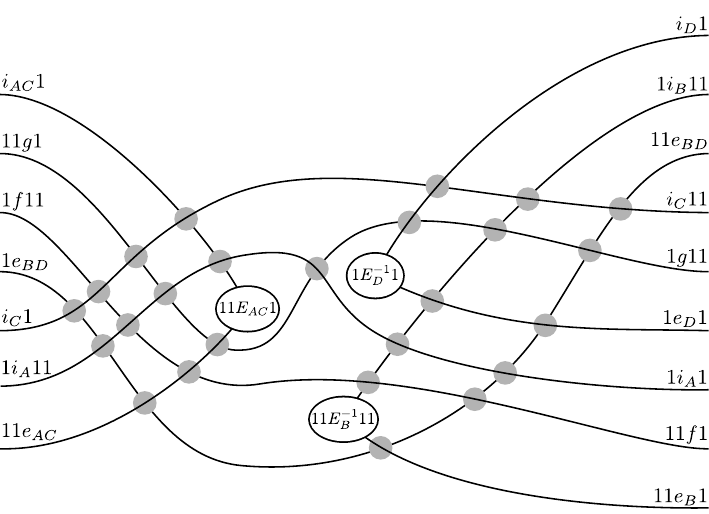}}

$$\begin{tabular}{@{}cc@{}}
\raisebox{0.45\adj}{$\chi_{f,g}:=$} & \includegraphics[width=90mm]{Pictures/Appendix/chinaturality.pdf}.
\end{tabular}$$

We now define three invertible modifications $\omega$, $\gamma$, and $\delta$. Given any objects $C$ and $D$ of $\mathfrak{C}$, we set $$\gamma_D:=F_D:\chi_{I,D}\Rightarrow Id_C,\ \delta_C:=F_C^{-1}:Id_C\Rightarrow \chi_{C,I}.$$ Finally, given any objects $B$, $C$, and $D$ in $\mathfrak{C}$, we define

\settoheight{\adj}{\includegraphics[width=60mm]{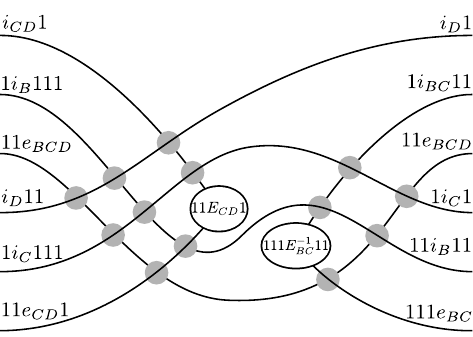}}

$$\begin{tabular}{@{}cc@{}}
\raisebox{0.45\adj}{$\omega_{B,C,D}:=$} & \includegraphics[width=60mm]{Pictures/Appendix/omega.pdf}.
\end{tabular}$$

\noindent It is easy to check that these invertible modifications satisfy the coherence conditions depicted in equations (HTA1) and (HTA2) of \cite{GPS}.
\end{proof}

\begin{Lemma}\label{lem:dual2natural}
Let $\mathfrak{C}$ and $\mathfrak{D}$ be monoidal 2-categories that have right duals, and let $F:\mathfrak{C}\rightarrow\mathfrak{D}$ be a monoidal 2-functor between them. There is a canonical monoidal 2-natural equivalence that witnesses the commutativity of the square $$\begin{tikzcd}[sep=small]
\mathfrak{C} \arrow[rr, "(-)^{\sharp}"] \arrow[dd, "F"'] &  & {\mathfrak{C}^{mop, 1op}} \arrow[dd, "{F^{mop,1op}}"] \\
 &  & \\
\mathfrak{D} \arrow[rr, "(-)^{\sharp}"'] &  &{\mathfrak{D}^{mop, 1op}}.
\end{tikzcd}$$
\end{Lemma}
\begin{proof}
For any $C$ in $\mathfrak{C}$, the equivalence in $\mathfrak{D}$ witnessing the commutativity of the above square is given by $$(F(C)^{\sharp}\Box F(e_C))\circ(F(C)^{\sharp}\Box \chi^F_{C,C^{\sharp}})\circ (i_{F(C)}\Box F(C^{\sharp})):F(C^{\sharp})\rightarrow F(C)^{\sharp}.$$ The remainder of the proof uses the same ideas as that of the previous proposition, so we leave the remaining details to the keen reader.
\end{proof}

\subsection{Diagrams for the proof of proposition \ref{prop:finitesemisimple}}\label{sub:FiniteSemisimpleProof}

\begin{figure}[!hbt]
    \centering %131.25
    \includegraphics[width=120mm]{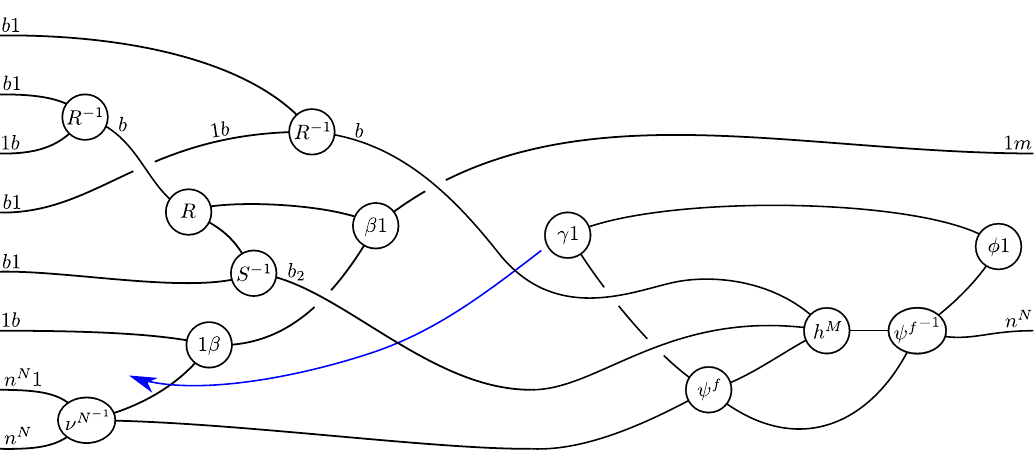}
    \caption{Cauchy Completeness (Part I.1)}
    \label{fig:cauchycompletenessI1}
\end{figure}

\begin{figure}[!hbt]
    \centering %112.5
    \includegraphics[width=108mm]{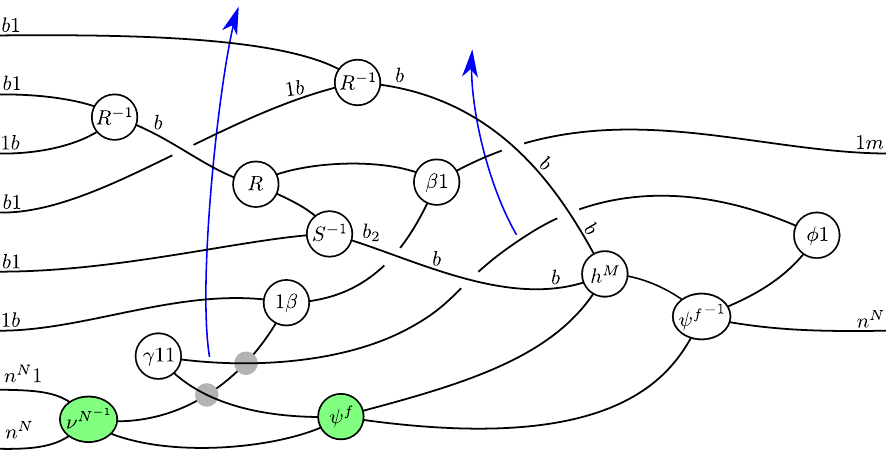}
    \caption{Cauchy Completeness (Part I.2)}
    \label{fig:cauchycompletenessI2}
\end{figure}

\begin{figure}[!hbt]
    \centering %135
    \includegraphics[width=120mm]{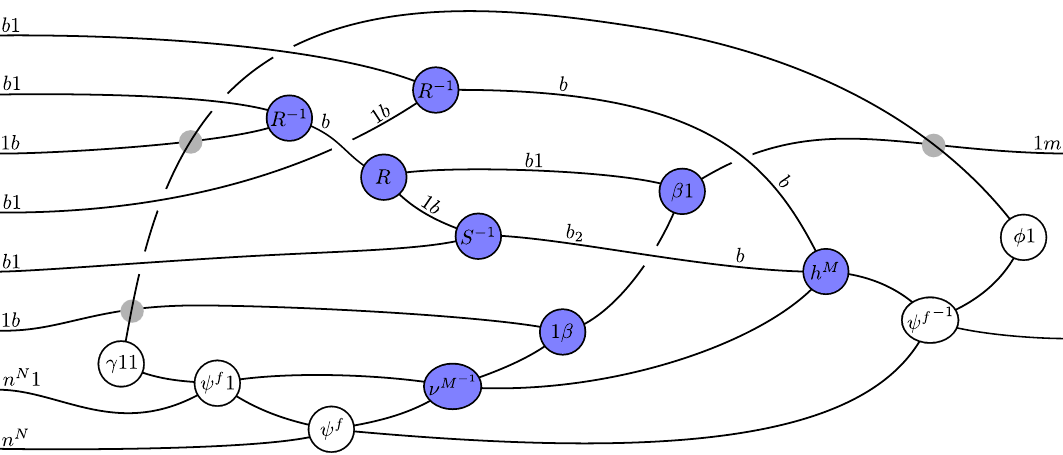}
    \caption{Cauchy Completeness (Part I.3)}
    \label{fig:cauchycompletenessI3}
\end{figure}

\begin{figure}[!hbt]
    \centering %112.5
    \includegraphics[width=110mm]{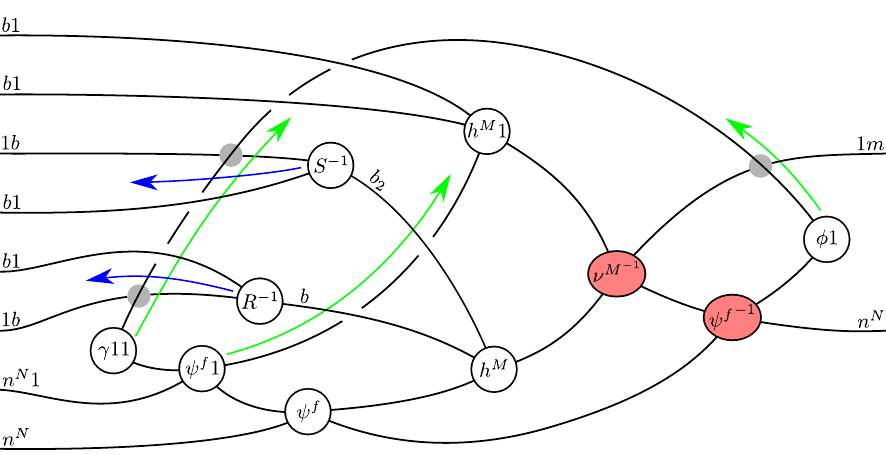}
    \caption{Cauchy Completeness (Part I.4)}
    \label{fig:cauchycompletenessI4}
\end{figure}

\begin{figure}[!hbt]
    \centering %120
    \includegraphics[width=110mm]{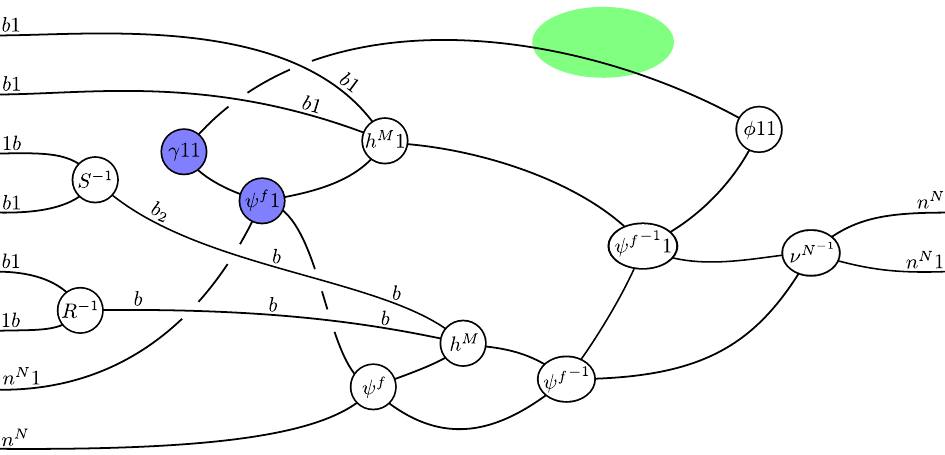}
    \caption{Cauchy Completeness (Part I.5)}
    \label{fig:cauchycompletenessI5}
\end{figure}

\begin{figure}[!hbt]
    \centering %127.5
    \includegraphics[width=115mm]{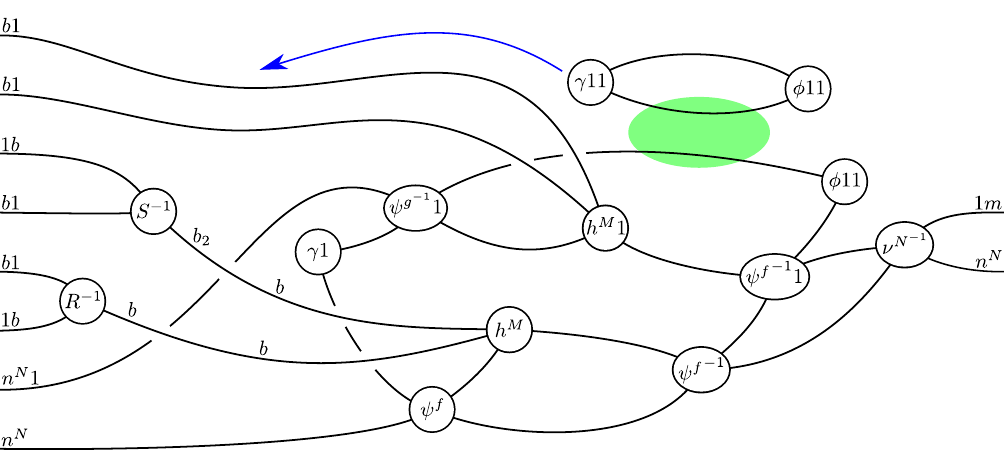}
    \caption{Cauchy Completeness (Part I.6)}
    \label{fig:cauchycompletenessI6}
\end{figure}

\begin{figure}[!hbt]
    \centering %127.5
    \includegraphics[width=115mm]{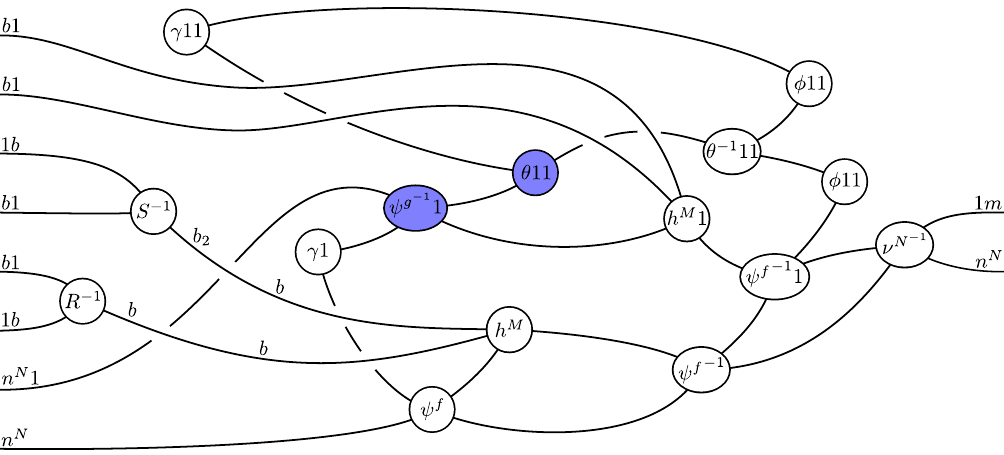}
    \caption{Cauchy Completeness (Part I.7)}
    \label{fig:cauchycompletenessI7}
\end{figure}

\begin{figure}[!hbt]
    \centering %127.5
    \includegraphics[width=115mm]{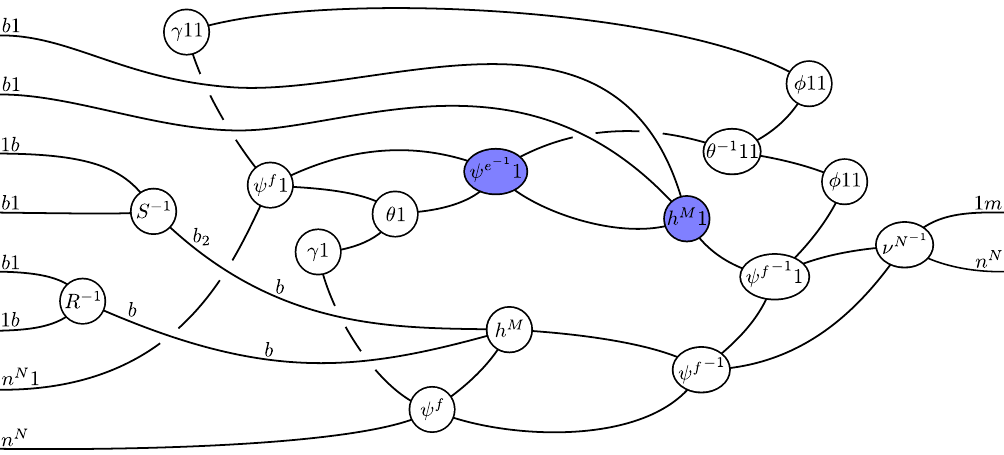}
    \caption{Cauchy Completeness (Part I.8)}
    \label{fig:cauchycompletenessI8}
\end{figure}

\begin{figure}[!hbt]
    \centering %120
    \includegraphics[width=112mm]{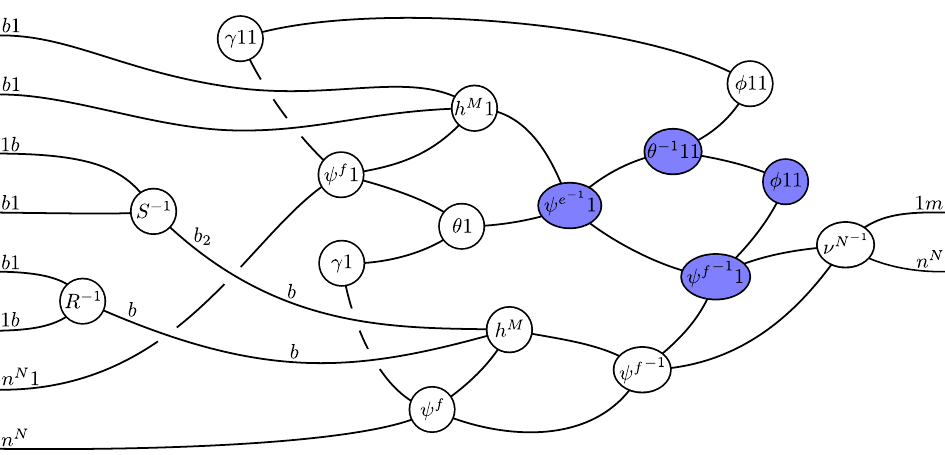}
    \caption{Cauchy Completeness (Part I.9)}
    \label{fig:cauchycompletenessI9}
\end{figure}

\begin{figure}[!hbt]
    \centering %120
    \includegraphics[width=112mm]{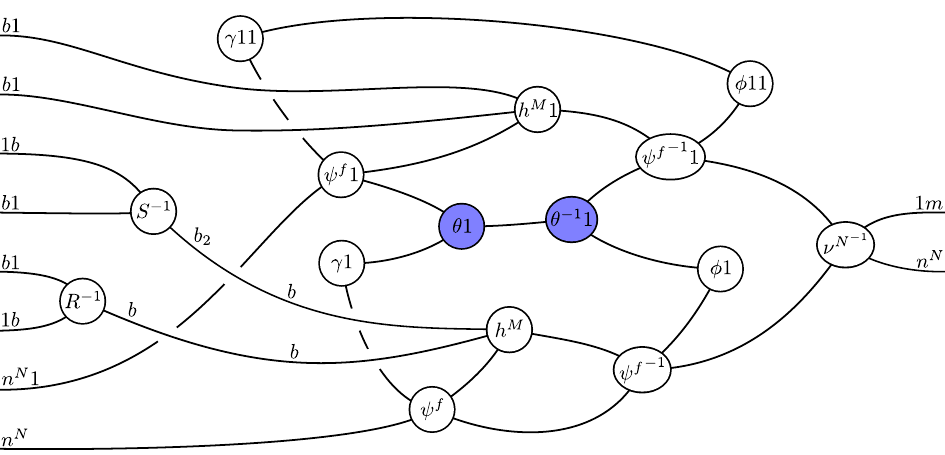}
    \caption{Cauchy Completeness (Part I.10)}
    \label{fig:cauchycompletenessI10}
\end{figure}

\begin{figure}[!hbt]
    \centering %120
    \includegraphics[width=112mm]{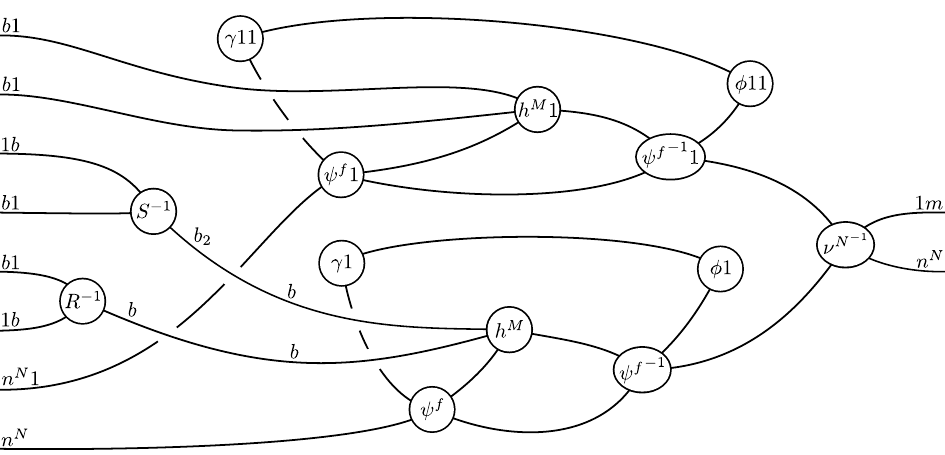}
    \caption{Cauchy Completeness (Part I.11)}
    \label{fig:cauchycompletenessI11}
\end{figure}

\FloatBarrier
    
     % Old
    
\bibliography{bibliography.bib}

\end{document}